\documentclass[preprint]{elsarticle}

\usepackage{amssymb}
\usepackage{latexsym}
\usepackage{amsthm}
\usepackage{enumerate}
\usepackage{epsfig}
\usepackage{graphicx}
\usepackage{color}
\usepackage{float}
\usepackage{subfigure}
\usepackage{amsmath}
\usepackage{makeidx}
\usepackage{fancyhdr}
\usepackage{lastpage}
\usepackage{url}

\usepackage{mathrsfs}

\usepackage{amsfonts}
\usepackage{dsfont}
\usepackage{mathtools}
\usepackage{bbm}
\usepackage{environ}
\usepackage{comment}
\usepackage{lineno,hyperref}
\modulolinenumbers[5]
\newcommand*{\normx}{\|{\cdot}\|}

\newcommand{\dm}{\mathrm{d}\mu}

\newcommand{\C}{\mathbb C}
\newcommand{\R}{\mathbb R}

\newcommand{\N}{\mathbb N}

\newcommand{\ei}{(e_i)_{i\in\N}}
\newcommand{\di}{(d_i)_{i\in\N}}

\newcommand{\fc}{\mathcal{F}C_b^{\infty}(}
\newcommand{\fce}{\mathcal{F}C_b^{\infty}}

\providecommand{\norm}[1]{\lVert#1\rVert} 
\providecommand{\abs}[1]{\lvert#1\rvert} 

\newtheorem{theorem}{Theorem}[section]

\newtheorem{hypo}[theorem]{Hypothesis}

\newtheorem{lemma}[theorem]{Lemma}

\newtheorem{corollary}[theorem]{Corollary}
\newtheorem{proposition}[theorem]{Proposition}

\newtheorem{definition}[theorem]{Definition}

\newtheorem{remark}[theorem]{Remark}

\hypersetup{pdfauthor={Name}}

\makeatletter
\def\ps@pprintTitle{%
  \let\@oddhead\@empty
  \let\@evenhead\@empty
  \def\@oddfoot{\reset@font\hfil\hfil\today}
  \let\@evenfoot\@oddfoot
}
\makeatother

\newcounter{case}

\DeclareUnicodeCharacter{2212}{-}

\journal{Journal of Functional Analysis}

\numberwithin{equation}{section}
\begin{document}
\begin{frontmatter}

\title{Essential m-dissipativity for generators of infinite-dimensional non-linear degenerate diffusion processes}

%
%

\author[coraut]{Benedikt Eisenhuth \fnref{myfootnote}}
\ead{eisenhuth@mathematik.uni-kl.de}
\author{Martin Grothaus \fnref{myfootnote}}
\ead{grothaus@mathematik.uni-kl.de}

\fntext[coraut]{Corresponding author.}
\fntext[myfootnote]{Department of Mathematics, TU Kaiserslautern, PO Box 3049, 67653 Kaiserslautern, Germany}

\begin{abstract}First essential m-dissipativity of an infinite-dimensional Ornstein-Uhlenbeck operator $N$, perturbed by the gradient of a potential, on a domain $\fce$ of finitely based, smooth and bounded functions, is shown. Our considerations allow unbounded diffusion operators as coefficients.
We derive corresponding second order regularity estimates for solutions $f$ of the Kolmogorov equation $\alpha f-Nf=g$, $\alpha \in (0,\infty)$, generalizing some results of Da Prato and Lunardi. Second we prove essential m-dissipativity for generators $(L_{\Phi},\fce)$ of infinite-dimensional non-linear degenerate diffusion processes. We emphasize that the essential m-dissipativity of $(L_{\Phi},\fce)$ is useful to apply general resolvent methods developed by Beznea, Boboc and Röckner, in order to construct martingale/weak solutions to infinite-dimensional non-linear degenerate diffusion equations. Furthermore, the essential m-dissipativity of $(L_{\Phi},\fce)$ and $(N,\fce)$, as well as the regularity estimates are essential to apply the general abstract Hilbert space hypocoercivity method from Dolbeault, Mouhot, Schmeiser and Grothaus, Stilgenbauer, respectively, to the corresponding diffusions.
\end{abstract}
%


\end{frontmatter}
\pagebreak
\section{Introduction}
The classical Langevin dynamics (compare \cite[Chapter 8.1]{physik1}) 
\begin{equation}
\begin{split}\label{endimLangevin}
	&\mathrm{d}X_t=Y_t\mathrm{d}t\\
	&\mathrm{d}Y_t=(-Y_t-D\Phi(X_t) )\mathrm{d}t+\sqrt{2}\mathrm{d}W_t,
\end{split}
\end{equation}
describes the evolution of the positions $X_t=(X_t^{(1)},...,X_t^{(n)})\in (\R^d)^n$ and velocities $Y_t=(Y_t^{(1)},...,Y_t^{(n)})\in (\R^d)^n$ of $n$ particles in dimension $d$, where $(W_t)_{t\geq 0}$ is a standard Brownian Motion in $(\R^d)^n$. I.e.~the velocity of the particles is subjected to friction and a stochastic perturbation. The $n$-particle potential $\Phi:(\R^d)^n\rightarrow\R$, with gradient $D\Phi$, affects the motion of the particles and can be used to model their interactions.

The equation has been studied under various aspects. In order to show exponential convergence to equilibrium of such type of non-coercive evolution equations, Cedric Villani developed the concepts of hypocoercivity, see \cite{vial}. Abstract hypocoercivity concepts for a quantitative descriptions of convergence rates are introduced in \cite{Dolbeaut}. These are translated to the corresponding stochastic equations, taking domain issues into account, in \cite{GrothausHypo}. In \cite{GrothausWangPoincare} these concepts have been further generalized to the case where only a weak Poincar\'e inequality is needed. In this case one obtains (sub-)exponential convergence rates.
%
%
Ergodicity and rate of convergence to equilibrium of the Langevin dynamics with singular potentials are elaborated e.g. in \cite{GrothausConrad} and \cite{david}. Recently, the dynamics and its hypocoercivity behavior is studied on abstract smooth manifolds, see \cite{LangevinManifold}. The latter articles cited above have in common that they study the associated Kolmogorov backward operator. Applying Itô's formula, the Kolmogorov backward operator associated to \eqref{endimLangevin}, also called Langevin operator, applied to $f\in C_0^{\infty}(\R^d\times\R^d)$ (w.l.o.g $n=1$) is given by
\begin{equation*}
L_{\Phi}f=\Delta_2f-\langle x,D_2f\rangle-\langle D\Phi,D_2f\rangle+\langle y,D_1f\rangle.
\end{equation*}
Here, $C_0^{\infty}(\R^d\times\R^d)$ denotes the space of compactly supported smooth (infinitely often differentiable) functions from $\R^d\times\R^d$ to $\R$, $x$ and $y$ the projection to the spatial and the velocity component, respectively, $\langle\cdot,\cdot\rangle$ the euclidean inner product, $\Delta_2$, $D_2$ the Laplacian and the gradient in the velocity component and $D_1$ the gradient in the spatial component. The key observation is the essential m-dissipativity of $(L_{\Phi},C_0^{\infty}(\R^d\times\R^d))$ defined in $L^2(\R^d\times\R^d,\mu^{\Phi},\R)$, where
\begin{align*}
\mu^{\Phi}=(2\pi)^{-\frac{d}{2}}e^{-\Phi(x)-\frac{1}{2}y^2}\lambda^d\otimes\lambda^d.
\end{align*}
We want to emphasis the degenerate structure of the Langevin equation, i.e.~the noise only appears in the velocity component. The degeneracy of the equation corresponds to the fact that the Laplacian in the definition of $(L_{\Phi},C_0^{\infty}(\R^d\times\R^d))$ is degenerate, i.e.~only acts in the velocity component. As the antisymmetric part of $(L_{\Phi},C_0^{\infty}(\R^d\times\R^d))$ in $L^2(\R^d\times\R^d,\mu^{\Phi},\R)$ contains first order differential operators in the spatial component and the symmetric part only differential operators in the velocity component, the operator $(L_{\Phi},C_0^{\infty}(\R^d\times\R^d))$ is non-sectorial.

\bigskip

This article addresses an infinite-dimensional generalization of the Langevin operator above.
In order to do that let $(U,(\cdot,\cdot)_U)$ and $(V,(\cdot,\cdot)_V)$ be two real separable Hilbert spaces. Consider the real separable Hilbert space $W=U\times V$ with inner product $(\cdot,\cdot)_W$ defined by
\begin{align*}
((u_1,v_1),(u_2,v_2))_W=(u_1,u_2)_U+(v_1,v_2)_V,\quad (u_1,v_1),(u_2,v_2)\in W. 
\end{align*}Denote by $\mathcal{B}(U)$ and $\mathcal{B}(V)$ the Borel $\sigma$-algebra on $U$ and $V$, on which we consider centered non-degenerate infinite-dimensional Gaussian measures $\mu_1$ and $\mu_2$, respectively. The measures are uniquely determined by their covariance operators $Q_1\in\mathcal{L}(U)$ and $Q_2\in\mathcal{L}(V)$. Furthermore we consider bounded linear operators $K_{12}\in\mathcal{L}(U;V)$, $K_{21}\in \mathcal{L}(U;V)$ and a symmetric bounded linear operator $K_{22}\in \mathcal{L}(V)$. For a given measurable potential $\Phi:U\rightarrow (-\infty,\infty]$, which is bounded from below, we set $\rho_{\Phi}=\frac{1}{c_{\Phi}}e^{-\Phi}$, where $c_{\Phi}=\int_U e^{-\Phi}\dm_1$ and consider the measure $\mu_1^{\Phi}=\rho_{\Phi}\mu_1$ on $(U,\mathcal{B}(U))$. 
On $(W,\mathcal{B}(W))$ we introduce the product measure
\begin{align*}
\mu^{\Phi}=\mu_1^{\Phi}\otimes\mu_2.
\end{align*}
The infinite-dimensional Langevin operator $(L_{\Phi},\fce)$ is defined by 
\begin{align*}
\fce  \ni f \mapsto L_{\Phi}f=S_{\Phi}f-A_{\Phi}f\in L^2(\mu^{\Phi}),
\end{align*} where for $f\in \fce$, $S_{\Phi}f$ and $A_{\Phi}f$ are given by
\begin{align*}
S_{\Phi}f=&\mathrm{tr}[K_{22}D^2_2f]-(v,Q_2^{-1}K_{22}D_2f)_V ,\\
A_{\Phi}f=&(u,Q_1^{-1}K_{21}D_2f)_U+(D \Phi(u),K_{21}D_2f)_U-(v,Q_2^{-1}K_{12}D_1f)_V.
\end{align*}Above, $\fce$ is a space of finitely based, smooth and bounded functions, see Definition \ref{spaces} and \ref{smooth}, below. Furthermore, $u$ and $v$ denotes the projections of $W$ to $U$ and $V$, respectively. One of the major challenges in this article is to show essential m-dissipativity of the infinite-dimensional Langevin operator.

We also address essential m-dissipativity and regularity estimates for infinite-dimensional Ornstein-Uhlenbeck operators, perturbed by the gradient of a potential $\Phi$. Indeed we fix a possible unbounded linear operator $(C,D(C))$ in $U$ and introduce the operator $(N,\fce)$ in $L^2(\mu_1^{\Phi})$ defined by
\begin{align*}
\fce \ni f \mapsto Nf=&\mathrm{tr}[CD^2f]-(u,Q_1^{-1}CD f)_U-(D\Phi,CDf)_U\in L^2(\mu_1^{\Phi}).
\end{align*}
In \cite{SobolevRegularity} and \cite{selfad} similar operators and corresponding regularity estimates are studied, but the results are restricted to bounded diffusion operators $(C,D(C))$ as coefficients.

The essential m-dissipativity of $(N,\fce)$ is useful in various applications. E.g.~to study stochastic quantization problems as in \cite[Section 4]{selfad} and to solve stochastic reaction diffusion equations as in \cite[Section 5]{SobolevRegularity}. In addition, essential m-dissipativity and related regularity estimates of such operators, will be essential for our planed application of the general abstract hypocoercivity method from \cite{GrothausHypo} to our infinite-dimensional setting. For this application it is needed to allow unbounded diffusion operators $(C,D(C))$ as coefficients in the definition of the perturbed Ornstein-Uhlenbeck operator $(N,\fce)$.


\bigskip
The organization of this article is as follows. First we fix notions and define several important spaces. Then properties of infinite-dimensional Gaussian measures are elaborated, especially the relation between finite and infinite-dimensional Gaussian measures in Lemma \ref{imagemeasure} and the integration by parts formula from Corollary \ref{IBP_formel_ad} are focused. In Theorem \ref{sobolev} we use the integration by parts formula to describe Sobolev spaces w.r.t.~infinite-dimensional Gaussian measures.

Section \ref{section_reg} introduces necessary conditions on $(C,D(C))$ (compare Hypothesis \ref{hypoCadv}) to obtain essential m-dissipativity of an Ornstein-Uhlenbeck operator with diffusion operator $(C,D(C))$ as coefficient. In Theorem \ref{essdispot} we perturb this Ornstein-Uhlenbeck operator by the gradient of a potential $\Phi:U\rightarrow (-\infty,\infty]$, which is in $W^{1,2}(U,\mu_1,\R)$ and bounded from below. If $D\Phi$ is strictly bounded by $\frac{1}{2\sqrt{\lambda_1}}$, where $\lambda_1$ is the biggest eigenvalue of $Q_1$, we obtain essential m-dissipativity of $(N,\fce)$. Note that the restriction to such potentials is due to the possible unboundedness of $(C,D(C))$. In the second part of this section we imitate the strategy used in \cite{SobolevRegularity} to derive an infinite-dimensional second order regularity estimate for $f\in \fce$ in terms of $g=\alpha f-Nf$, $\alpha\in (0,\infty)$.

In Section \ref{section_lan}, we deal with the essential m-dissipativity of $(L_{\Phi},\fce)$. First we consider the case where $\Phi=0$. We decompose our infinite-dimensional Langevin operator into countable finite-dimensional ones, to use perturbation arguments for finite-dimensional Langevin operators as described in \cite{GrothausConrad} and \cite{DissConrad}. We also derive regularity estimates similar to the ones from Section \ref{section_reg}. Finally we consider potentials $\Phi$ as in Theorem \ref{ess_diss_Lphi} and use the Neumann-Series theorem to obtain essential m-dissipativity of $(L_{\Phi},\fce)$ in $L^2(\mu^{\Phi})$. During the whole section we assume Hypothesis \ref{stass}, which is the key to the decomposition described above.

Applications of the results, we derived in Section \ref{section_reg} and Section \ref{section_lan}, are discussed in the last section. We propose an infinite-dimensional non-linear degenerate stochastic differential equation, see \eqref{sdeex}. With the results we achieved in this article and the resolvent methods from \cite{boboc} we plan to solve it. Moreover we elaborate, how the essential m-dissipativity of $(N,\fce)$ can be used to show hypocoercivity of the semigroup $(T_t)_{t\geq 0}$ generated by $(L_{\Phi},D(L_{\Phi}))$ and how hypocoercivity of $(T_t)_{t\geq 0}$ is related to the long time behavior of the process solving \eqref{sdeex}.
The main results obtained in this article are summarized in the following list:

\begin{enumerate}
\item[$\bullet$]We prove essential m-dissipativity of perturbed Ornstein-Uhlenbeck operators $(N,\fce)$ in $L^2(\mu_1^{\Phi})$, with possible unbounded diffusions $(C,D(C))$ as coefficients, see Theorem \ref{essdispot}. There Hypothesis \ref{hypoCadv} is assumed and perturbations by the gradient of a potential $\Phi\in W^{1,2}(U,\mu_1,\R)$, which is bounded from below, are considered. In addition, an appropriate bound for the gradient of $\Phi$, i.e.~$\norm{D\Phi}_{L^{\infty}(\mu_1)}<\frac{1}{2\sqrt{\lambda_1}}$, where $\lambda_1$ is the biggest eigenvalue of $Q_1$ (see Theorem \ref{essdispot}), is needed.
\item[$\bullet$]Considering potentials $\Phi\in W^{1,2}(U,\mu_1,\R)$, which are convex, bounded from below, lower semicontinuous and with $\int_U \norm{D\Phi}_{U}^p\dm_1< \infty$ for some $p\in (2,\infty)$ as in Hypothesis \ref{hypo_potential_convex}, we provide second order regularity estimate for $f\in \fce$ in terms of $g=\alpha f-Nf$, $\alpha\in (0,\infty)$. Indeed by Theorem \ref{theorem_reg} it holds
\begin{align*}
\int_U \mathrm{tr}[(CD^2f)^2]+\norm{Q_1^{-\frac{1}{2}}CDf}_U^2\dm_1^{\Phi}\leq 4\int_U g^2\dm_1^{\Phi},
\end{align*}where $(C,D(C))$ is the possible unbounded diffusion coefficient in the definition of $(N,\fce)$.
\item[$\bullet$]Essential m-dissipativity of the infinite-dimensional Langevin operator $(L_{\Phi},\fce)$ in $L^2(\mu^{\Phi})$ is shown in Theorem \ref{ess_diss_Lphi}. We consider potentials $\Phi$ in $W^{1,2}(U,\mu_1,\R)$, which are bounded from below, with $\norm{D\Phi}_{L^{\infty}(\mu_1)}<\frac{1}{2\sqrt{\nu_1}}$, where $\nu_1$ is the biggest eigenvalue of $Q_2$ and assume Hypothesis \ref{stass}.
\end{enumerate}

\section{Notations and preliminaries}
Let $U$ and $V$ be two real separable Hilbert spaces with inner products $(\cdot,\cdot)_U$ and $(\cdot,\cdot)_V$, respectively. The induced norms are denoted by $\normx_U$ and $\normx_V$. 
The set of all linear bounded operators from $U$ to $U$ and from $U$ to $V$ are denoted by $\mathcal{L}(U)$ and $\mathcal{L}(U;V)$. By $\mathcal{L}^+(U)$ we shall denote the subset of $\mathcal{L}(U)$ consisting of all nonnegative symmetric operators. The subset of operators in $\mathcal{L}^+(U)$ of trace class is denoted by $\mathcal{L}^+_1(U)$ and the set of Hilbert-Schmidt operators by $\mathcal{L}_2(U)$.

Suppose we have $J\in \mathcal{L}^+(U)$. If $J$ is injective it is reasonable to talk about the the inverse of $J:U\rightarrow J(U)$, which will be denoted by $J^{-1}$. Due to \cite[Proposition 4.4.8.]{Neumann} there exists a unique operator $J^{\frac{1}{2}}\in \mathcal{L}^+(U)$ such that $(J^\frac{1}{2})^2=J$. If $J^{-1}$ exists, so does $(J^{\frac{1}{2}})^{-1}$, in this case we denote $(J^{\frac{1}{2}})^{-1}$ by $J^{-\frac{1}{2}}$.
By $\mathcal{B}(U)$ we denote the Borel $\sigma$-algebra, i.e.~the $\sigma$-algebra generated by the open sets in $(U,(\cdot,\cdot)_U)$. The euclidean inner product and induced norm is denoted by $\langle\cdot,\cdot\rangle$ and $\abs{\cdot}$, respectively.

For a given measure space $(\Omega,\mathcal{A},m)$ and a Banach space $Y$ we denote by $L^p(\Omega,m,Y)$, $p\in [0,\infty]$ the Hilbert space of equivalence classes of $\mathcal{A}$-$\mathcal{B}(Y)$ measurable and $p$-integrable functions. The corresponding norm is denoted by $\normx_{L^p(\Omega,m,Y)}$. If $p=2$, the norm is induced by an inner product denoted by $(\cdot,\cdot)_{L^2(\Omega,m,Y)}$. In case ($\Omega,\mathcal{A})$ is clear from the context and $Y=\R^n$ for some $n\in \N$, we also write $L^2(m)$ instead of $L^2(\Omega,m,\R^n)$. By $\lambda^n$, $n\in\N$, we denote the Lebesgue measure on $(\R^n,\mathcal{B}(\R^n))$.

On the measurable space $(U,\mathcal{B}(U))$ we consider an infinite-dimensional non-degenerate Gaussian measure $\mu_1$ with covariance operator $Q_1\in\mathcal{L}^+_1(U)$. Since the measure is non-degenerate the operator $Q_1$ is injective and therefore positive. For the definition and construction of these measures we refer to the first chapter of \cite{DaPrato_intro_infinite_dim_analysis}.

In the next lemma we discuss the important relation between finite and infinite-dimensional Gaussian measures. A proof can be found in \cite[Corollary 1.19]{DaPrato_intro_infinite_dim_analysis}.
\begin{lemma}\label{imagemeasure}Given $n\in\N$ and elements $l_1,...,l_n\in X$. The image measure $\mu_1^n$ of $\mu_1$ under the map
\begin{align*}
U \ni u \mapsto \big((l_1,u)_U,...(l_n,u)_U \big)\in\R^n
\end{align*}is the centered Gaussian measure  on $(\R^n,\mathcal{B}(\R^n))$ with covariance matrix $Q_{1,n}=((Q_1l_i,l_j)_U)_{ij=1,...,n}$.\\
If $l_1,...,l_n$ is an orthonormal system of eigenvectors of $Q_1$ with corresponding eigenvalues $\lambda_1,...,\lambda_n$, the covariance matrix $Q_{1,n}$ of $\mu_1^n$ is given by the diagonal matrix $\mathrm{diag}(\lambda_1,...,\lambda_n)$.
\end{lemma}

During this article we have to perform explicit calculations of integrals with respect to Gaussian measures including monomials of order 2 and 4, therefore the following lemma is useful. To proof it, apply Lemma \ref{imagemeasure} and Isserlis formula from \cite{isserlis}.
\begin{lemma}\label{Gaussianmoments}For $l_1,l_2,l_3,l_4\in U$ it holds
\begin{align*}
&\int_U (u,l_1)_U(u,l_2)_U\mathrm{d}\mu_1(u)=(Q_1l_1,l_2)_U\quad \text{and}\\
&\int_U (u,l_1)_U(u,l_2)_U(u,l_3)_U(u,l_4)_U\mathrm{d}\mu_1(u)\\
&=(Q_1l_1,l_2)_U(Q_1l_3,l_4)_U+(Q_1l_1,l_3)_U(Q_1l_2,l_4)_U+(Q_1l_1,l_4)_U(Q_1l_2,l_3)_U.
\end{align*}
\end{lemma}
To cover more general situations we consider a measurable potential $\Phi:U\mapsto (-\infty,\infty]$, which is bounded from below. During the paper we will assume more or less restrictive assumptions on the potential. 
As in the introduction we set $\rho_{\Phi} =\frac{1}{c_{\Phi}}e^{-\Phi}$, where $c_{\Phi}=\int_U e^{-\Phi}\mathrm{d}\mu_1$. On $(U,\mathcal{B}(U))$ we consider the measure $\mu_1^{\Phi}$ defined by
 \begin{align*}
 \mu_1^{\Phi}=\rho_{\Phi}\mathrm{d}\mu_1.
 \end{align*}I.e.~a measures having a density with respect to the infinite-dimensional Gaussian measure $\mu_1$. We fix an orthonormal basis $B_U=\di$ of $U$.
\begin{definition}\label{spaces}For $n\in\N$, set $B_U^n=\mathrm{span}\lbrace d_1,...,d_n\rbrace$. The orthogonal projection from $U$ to $B^n_U$ is denoted by $\overline{P}_n$ and the corresponding coordinate map by $P_n$, i.e. we have for all $u\in U$
\begin{align*}
\overline{P}_n(u)=\sum_{i=1}^n (u,d_i)_Ud_i\quad \text{and}\quad P_n(u)=\big( (u,d_1)_U,...,(u,d_n)_U\big).
\end{align*}
Let $C_b^{\infty}(\R^n)$ be the space of all bounded smooth (infinitely often differentiable) real-valued functions on $\R^n$. The space of finitely based smooth and bounded functions on $U$, is defined by
\begin{align*}
\fc B_U )&=\lbrace U\ni u\mapsto \varphi(P_m(u))\in \R\mid m\in\N, \; \varphi\in C_b^{\infty}(\R^{m})\rbrace.
\end{align*}The subset of functions only depending on $n$-directions is defined correspondingly by
\begin{align*}
\fc B_U,n)&=\lbrace U\ni u\mapsto \varphi(P_n(u))\in \R\mid  \varphi\in C_b^{\infty}(\R^{n})\rbrace.
\end{align*}For later use we define 
\begin{align*}
L^2_{B_{U_n}}(U,\mu_1,\R)=\lbrace U\ni u\mapsto f(P_n(u))\in\R\mid f\in L^2(\R^n,\mu_1^n,\R)\rbrace,
\end{align*}where $\mu_1^n$ is the image measure of $\mu_1$ under $P_n$. With the inner product induced by $L^2(U,\mu_1,\R)$ this space is a Hilbert space. 
\end{definition}
A very useful density result, proved in \cite[Lemma 2.2]{SobolevRegularity}, is stated in the next lemma.

\begin{lemma}\label{density}The function spaces $\fc B_U)$ and $\fc B_U,n)$ are dense in $L^2(U,\mu_1,\R)$ and $L^2_{B_{U_n}}(U,\mu_1,\R)$, respectively.
\end{lemma}

\begin{remark}\label{derivatives}
Given a Frech\'{e}t differentiable function $f:U\rightarrow \R$. For $u\in U$
we denote by $Df(u)\in U$ the gradient of $f$ in $u$.
Analogously for a two times Frech\'{e}t differentiable function $f:U\rightarrow \R$, we identify $D^2f(u)\in L(U)$ with the second order Frech\'{e}t derivative in $u\in U$.
For $i,j\in \N$ we denote by $\partial_if(u)=(Df(u),d_i)_U$ the partial derivative in the direction of $d_i$ and by $\partial_{ij}f(u)=(D^2f(u)d_i,d_j)_U$ the second order partial derivative in the direction of $d_i$ and $d_j$ .
\end{remark}

We continue this section with the important integration by parts formula for infinite-dimensional Gaussian measures, with and without densities. We assume that $B_U=\di$ is an orthonormal basis of eigenvectors of $Q_1$ with corresponding eigenvalues $(\lambda_i)_{i\in\N}\subset (0,\infty)$. W.l.o.g. we assume that $(\lambda_i)_{i\in\N}$ is decreasing to zero.\\
Since $Q_1$ is injective, the inverse $Q_1^{-1}$ of $Q_1:U\rightarrow Q_1(U)$ exists. Obviously it holds
\begin{align*}
Q_1^{-1}d_i=\frac{1}{\lambda_i}d_i,\quad i\in\N,
\end{align*}and therefore it is reasonable to define the operator $Q_1^{-\frac{1}{2}}$ on $\bigcup_{n\in\N}B_U^n$ uniquely characterized by 
\begin{align*}
Q_1^{-\frac{1}{2}}d_i=\frac{1}{\sqrt{\lambda_i}}d_i,\quad i\in\N.
\end{align*}
\begin{theorem}\label{IBP_formel}
For $f,g\in \fc B_U)$ and $i\in\N$, it holds the integration by parts formula
\begin{align}
\int_U \partial_ifg\mathrm{d}\mu_1=-\int_U f\partial_ig\mathrm{d}\mu_1+\int_U (u,Q^{-1}d_i)_U fg\mathrm{d}\mu_1.
\end{align}
\begin{proof}
Apply Lemma \ref{imagemeasure}, use the standard integration by parts formula and use Lemma \ref{imagemeasure} again. For a more detailed proof see \cite[Lemma 10.1]{DaPrato_intro_infinite_dim_analysis}.
\end{proof}
\end{theorem}
Fix a possible unbounded linear operator $(C,D(C))$ on $U$ fulfilling the following hypothesis.
\begin{hypo}\label{hypo_C_one}$(C,D(C))$ is nonnegative symmetric and for all $n\in\N$, it holds
\begin{align*}
B_U^n\subset D(C)\quad C(B_U^n)\subset B_U^n.
\end{align*}
\end{hypo}
\begin{remark}\label{domain}Note that for given $n\in \N$ and $f=\varphi(P_n(\cdot))\in\fc B_U)$ one has $Df=\sum_{i=1}^n\partial_i\varphi(P_n(\cdot))d_i\in B_U^n$. Hence the hypothesis above ensure that expressions like $CDf$, $Q_1^{-\frac{1}{2}}CDf$ or $Q_1^{-1}CDf$ are well-defined.
\end{remark}

\begin{theorem}\label{sobolev}Assume that Hypothesis \ref{hypo_C_one} hold.
The operators 
\begin{align*}
D:\fc B_U)&\rightarrow L^2(U,\mu_1,U)\\
CD:\fc B_U)&\rightarrow L^2(U,\mu_1,U)\\
Q_1^{-\frac{1}{2}}CD:\fc B_U)&\rightarrow L^2(U,\mu_1,U)\\
(CD,CD^2):\fc B_U)&\rightarrow L^2(U,\mu_1,U)\times L^2(U,\mu_1,\mathcal{L}_2(U))
\end{align*}
are closable in $L^2(U,\mu_1,\R)$. 
\begin{proof}As the proof of closability for the first three operators are essentially the same, we restrict ourself to the third and the last. 
Let $(f_n)_{n\in\N}\subset \fc B_U)$ converge to $0$ in $L^2(U,\mu_1,\R)$ and be such that $Q_1^{-\frac{1}{2}}CDf_n\rightarrow F$ in $L^2(U,\mu_1,U)$ as $n\rightarrow\infty$. For $k\in \N$ we have by the invariance properties of $(C,D(C))$ and the fact that $\di$ is an orthonormal basis of eigenvectors of $Q_1$
\begin{align*}
(Q_1^{-\frac{1}{2}}CDf_n,d_k)_U=\sum_{l=1}^{\infty}(Q_1^{-\frac{1}{2}}Cd_l,d_k)_U\partial_lf_n=\sum_{l=1}^{k}(Q_1^{-\frac{1}{2}}Cd_l,d_k)_U\partial_lf_n.
\end{align*}
For an arbitrary $g\in \fc B_U)$ we obtain by the integration by parts formula
\begin{align*}
&\int_U (Q_1^{-\frac{1}{2}}CDf_n,d_k)_U g \mathrm{d}\mu_1\\
=&-\sum_{l=1}^{k}(Q_1^{-\frac{1}{2}}Cd_l,d_k)_U\int_U f_n(\partial_lg-(u,Q_1^{-1}d_l)_Ug)\mathrm{d}\mu_1.
\end{align*}Observe that $g$ and $\partial_lg-(u,Q_1^{-1}d_l)_Ug$ are in $L^2(U,\mu,\R)$ and therefore
\begin{align*}
\int_X (F,d_k)_U g \mathrm{d}\mu_1=0.
\end{align*}By the density of $\fc B_U)$ in $L^2(U,\mu_1,\R)$ we conclude $(F,d_k)_U=0$ for all $k\in\N$, hence finally $F=0$.

To show that the fourth operator is closable we proceed similarly. Indeed, let $(f_n)_{n\in\N}\subset \fc B_U)$ converge to $0$ in $L^2(U,\mu_1,\R)$ and be such that $CDf_n\rightarrow F$ in $L^2(U,\mu_1,U)$ and $CD^2f_n\rightarrow A$ in $L^2(U,\mu_1,\mathcal{L}_2(U))$, as $n\rightarrow\infty$. As above $F=0$. Now for $k,l\in\N$, we have
\begin{align*}
(CD^2f_nd_k,d_l)_U=\sum_{i=1}^l (Cd_l,d_i)_U\partial_{ki}f_n.
\end{align*}
Hence for arbitrary $g\in \fc B_U)$, we obtain by the integration by parts formula
\begin{align*}
\int_U(CD^2f_nd_k,d_l)_Ug\mathrm{d}\mu_1&=\sum_{i=1}^l (Cd_l,d_i)_U\int_U\partial_{ki}f_ng\mathrm{d}\mu_1\\
&=-\sum_{i=1}^l (Cd_l,d_i)_U\int_U \partial_if_n(\partial_kg-(u,Q_1^{-1}d_k)_Ug)\mathrm{d}\mu_1\\
&=-\int_U (d_l,CDf_n)_U(\partial_kg-(u,Q_1^{-1}d_k)_Ug)\mathrm{d}\mu_1.
\end{align*}Arguing as in the first part we observe $(Ad_k,d_l)_U=0$ in $L^2(U,\mu_1,\R)$, implying $A=0$ in $L^2(U,\mu_1,\mathcal{L}_2(U))$.
\end{proof}
\end{theorem}
By Theorem \ref{sobolev} it is reasonable to define $W^{1,2}(U,\mu_1,\R)$, $W^{1,2}_{C}(U,\mu_1,\R)$, $W^{1,2}_{Q_1^{-\frac{1}{2}}C}(U,\mu_1,\R)$ and $W^{2,2}_{C}(U,\mu_1,\R)$ as the domain of the closures of $D$, $CD$, $Q_1^{-\frac{1}{2}}CD$ and $(CD,CD^2)$ in $L^2(U,\mu_1,\R)$, respectively. We still denote the closures of the differential operators from the theorem above by $D$, $CD$, $Q_1^{-\frac{1}{2}}CD$ and $(CD,CD^2)$. By \cite[Proposition 10.6]{DaPrato_intro_infinite_dim_analysis} every bounded function $f:U\rightarrow \R$ with bounded Frech\'{e}t derivative is in $W^{1,2}(U,\mu_1,\R)$ and the classical gradient of $f$ coincides with $Df$ in $L^2(U,\mu_1,\R)$.

\begin{remark}\label{IBP_formel_ad}Adapting the proof of \cite[Lemma 9.2.5]{Secondorder} one can show that the integration by parts formula from Theorem \ref{IBP_formel} also holds in the case that $f,g\in W^{1,2}(U,\mu_1,\R)$. 
\end{remark}
Invoking the remark above the following integration by parts formula for the measure $\mu^{\Phi}$ is valid.
\begin{corollary}\label{IBP_potential}Assume the potential $\Phi:U\rightarrow (-\infty,\infty]$ is bounded from below and in $W^{1,2}(U,\mu_1,\R)$.
For $f,g\in \fc B_U)$ and $i\in\N$, it holds the integration by parts formula
\begin{align*}
\int_U \partial_ifg\mathrm{d}\mu_1^{\Phi}=&-\int_U f\partial_ig\mathrm{d}\mu_1^{\Phi}+\int_U (u,Q^{-1}d_i)_U f(u)g(u)\mathrm{d}\mu_1^{\Phi}+\int_U \partial_i\Phi fg\mathrm{d}\mu_1^{\Phi}.
\end{align*}
\end{corollary}

\section{Perturbed Ornstein-Uhlenbeck operators and corresponding regularity estimates}\label{section_reg}
This section is devoted to an infinite-dimensional Ornstein-Uhlenbeck operator, perturbed by the gradient of the potential $\Phi$. As already mentioned in the introduction, such operators naturally occur during the application of the abstract Hilbert space hypocoercivity method. Also an infinite-dimensional regularity estimate is derived in the second part of this section. The proof of such estimates is motivated by the results from \cite{SobolevRegularity}, where Giuseppe Da Prato and Alessandra Lunardi investigated Sobolev regularity for a class of second order elliptic partial differential equations in infinite-dimensions. As before $\mu_1$ is a centered non-degenerate Gaussian measure on the real separable Hilbert space $U$ with covariance operator $Q_1\in \mathcal{L}_1^+(U)$. We fix an orthonormal basis of eigenvectors $B_U=\di$ of $Q_1$. W.l.o.g. the corresponding sequence of eigenvalues $(\lambda_i)_{i\in\N}$ decreases to zero. We start with the operator $(N_0,\fc B_U))$, defined in $L^2(U,\mu_1,\R)$ by
\begin{align*}
\fc B_U)\ni f\mapsto N_0f=\mathrm{tr}[CD^2f]-(u,Q^{-1}CD f)_U\in L^2(U,\mu_1,\R),
\end{align*}where we assume that $(C,D(C))$ is a possible unbounded linear operator on $U$. Since we allow such unbounded diffusions as coefficients, we cannot use general results from \cite{SobolevRegularity}  or \cite[Section 10]{Secondorder}. Assuming Hypothesis \ref{hypoCadv} below, ensures that the expressions $\mathrm{tr}[CD^2f]$ and $Q_1^{-1}CD f$ are reasonable.

At this point we have to mention that the operator $N_0$ is well-defined in the sense that two representatives of the same equivalence class yield the same output. To see this, note that the measure $\mu_1$ has full topological support, i.e. the smallest closed measurable set with full measure is $U$. The proof of this statement can be found in \cite{support}, it relies on the fact that we assumed that the Hilbert space $U$ is separable.

During this section we permanently assume the following hypothesis.
\begin{hypo}\label{hypoCadv}$ $
\begin{enumerate}
\item $(C,D(C))$ is symmetric and nonnegative. 
\item For all $n\in\N$ it holds
\begin{align*}
B_U^n\subset D(C)\quad C(B_U^n)\subset B_U^n.
\end{align*}
\item The operator $(-Q_1^{-1}C,D(Q_1^{-1}C))$ on $U$ is of negative type, i.e. 
there is a $\tau_1 > 0$ such that for all $u\in D(Q_1^{-1}C)=\lbrace u\in D(C)\vert\; Cu\in D(Q_1^{-1})\rbrace$
\begin{align*}
(-Q_1^{-1}Cu,u)_U\leq -\tau_1\norm{u}^2.
\end{align*}
\end{enumerate} 
\end{hypo}

\begin{theorem}\label{ornstein_ess}The operator $(N_0,\fc B_U))$ is 
\begin{enumerate}
\item dissipative in $L^2(U,\mu_1,\R)$, with $(N_0f,g)_{L^2(U,\mu_1,\R)}=\int_U -(CDf,Dg)_U\dm_1$ for all $f,g\in \fc B_U)$ and fulfills
\item the dense range condition $\overline{(Id-N_0)(\fc B_U))}=L^2(U,\mu_1,\R)$,
\end{enumerate}i.e.~is essentially m-dissipative in $L^2(U,\mu_1,\R)$. The resolvent in $\alpha\in (0,\infty)$ of the closure $(N_0,D(N_0))$ is denoted by $R(\alpha,N_0)$.
\begin{proof}
The first item of the statement follows by the integration by parts formula from Theorem \ref{IBP_formel} together with the invariance properties of $(C,D(C))$. For the second statement we fix $n\in \N$ and set
\begin{align*}
C_n=((Cd_i,d_j)_U)_{ij=1}^n\quad\text{and}\quad Q_{1,n}=((Q_1d_i,d_j)_U)_{ij=1}^n=\mathrm{diag}(\lambda_1,...,\lambda_n).
\end{align*}For $f=\varphi(P_n(\cdot))\in\fc B_U,n)$ and $u\in U$ we get by the invariance properties of $(C,D(C))$
\begin{align*}
N_0f(u)=\mathrm{tr}[C_nD^2\varphi(P_n(u))]-\langle P_n(u),Q_{1,n}^{-1}C_nD\varphi(P_n(u))\rangle.
\end{align*}
It is therefore natural to consider the operator $(N_{0,n},C_b^{\infty}(\R^n))$ defined by
\begin{align*}
C_b^{\infty}(\R^n)\ni \varphi \mapsto N_{0,n}\varphi=\mathrm{tr}[C_nD^2\varphi]-\langle\cdot,Q^{-1}_nC_nD\varphi\rangle\in L^2(\R^n,\mu_1^n,\R)
\end{align*} 
It is well known (compare \cite[Proposition 10.2.1]{Secondorder}) 
that $(N_{0,n},C_b^{\infty}(\R^n))$ is essential m-dissipative in $L^2(\R^n,\mu_1^n,\R)$, hence $(Id-N_{0,n})(C_b^{\infty}(\R^n))$ is dense in $L^2(\R^n,\mu_1^n,\R)$. 

Given $\varepsilon >0$ and $h=g(P_n(\cdot))\in L^2_{B_{U_n}}(U,\mu_1,\R)$. As $g$ is in $L^2(\R^n,\mu_1^n,\R)$ and $(Id-N_{0,n})(C_b^{\infty}(\R^n))$ is dense in $L^2(\R^n,\mu_1^n,\R)$ we find a $\varphi\in C_b^{\infty}(\R^n)$ such that
\begin{align*}
\norm{(Id-N_{0,n})\varphi - g}_{L^2(\R^n,\mu_1^n,\R)}<\varepsilon.
\end{align*}
Using Lemma \ref{imagemeasure} we obtain
\begin{align*}
\norm{(Id-N_0)\varphi(P_n(\cdot))-h}_{L^2_{B_{U_n}}(U,\mu_1,\R)}&=\norm{(Id-N_0)\varphi(P_n(\cdot))-h}_{L^2(U,\mu_1,\R)}\\
&=\norm{(Id-N_{0,n})\varphi - g}_{L^2(\R^n,\mu_1^n,\R)}<\varepsilon.
\end{align*}
In other words $(Id-N_0)(\fc B_U,n))$ is dense in $L^2_{B_{U_n}}(U,\mu_1,\R)$.

Finally we use the result above to show that $(Id-N_0)(\fc B_U))$ is dense in $L^2(U,\mu_1,\R)$. Indeed let $\varepsilon>0$ and take an element $h\in L^2(U,\mu_1,\R)$. By Lemma \ref{density} we find a $g\in \fc B_U,n)\subset L^2_{B_{U_n}}(U,\mu_1,\R)$ such that
\begin{align*}
\norm{h-g}_{L^2(U,\mu_1,\R)}< \frac{\varepsilon}{2}.
\end{align*}As $\overline{(Id-N_0)(\fc B_U,n))}=L^2_{B_{U_n}}(U,\mu_1,\R)$ we find a $f\in \fc B_U,n)$ such that 
\begin{align*}
\norm{(Id-N_0)f-g}_{L^2(U,\mu_1,\R)}< \frac{\varepsilon}{2}.
\end{align*}Hence the triangle inequality yields
\begin{align*}
\norm{(Id-N_0)f-h}_{L^2(U,\mu_1,\R)}< \varepsilon.
\end{align*}
Invoking the famous Lumer-Phillips theorem we obtain that $(N_0,\fc B_U))$ is essential m-dissipative in $L^2(U,\mu_1,\R)$.
\end{proof}
\end{theorem}
Before we go ahead and perturb $(N_0,\fc B_U))$ we need a $L^2(\mu_1)$ regularity estimate for the first and second order derivatives of a function $f\in\fc B_U)$ in terms of $g\in L^2(U,\mu_1,\R)$, where
\begin{equation}\label{equationUreg}
\alpha f-N_0f=g,
\end{equation}for a given $\alpha \in (0,\infty)$. 
\begin{theorem}\label{regestimates}
Suppose we have $f\in\fc B_U)$ and $g=\alpha f-N_0f$, $\alpha\in (0,\infty)$, as in Equation \eqref{equationUreg} above. It holds $g\in W^{1,2}(U,\mu_1,\R)$ and the identities 
\begin{align*}
&\int_U\alpha f^2+(CD f,D f)_U\dm_1=\int_Ugf\dm_1\\
 &\int_U\alpha \norm{CDf}^2_U+\norm{Q_1^{-\frac{1}{2}}CDf}^2_U+\mathrm{tr}[(CD^2f)^2]\dm_1
 =\int_U (D g,CD f)_U\dm_1,
\end{align*}are valid. In particular it holds
\begin{align*}
 \int_U\norm{Q_1^{-\frac{1}{2}}CDf}^2_U+\mathrm{tr}[(CD^2f)^2]\dm_1=\int_U (N_0f)^2\dm_1.
\end{align*}
\begin{proof}Due to the definition of $N_0$ and the fact that $f$ is in $\fc B_U)$ it is easy to see that $g$ is infinitely often differentiable. As $Dg$ is in $L^2(U,\mu_1,U)$ and has at most linear growth (compare equation \eqref{partialg}), an approximation argument shows that $g$ is in $W^{1,2}(U,\mu_1,\R)$.\\
To show the first equation, we multiply \eqref{equationUreg} with $f$ and integrate over $U$ with respect to $\mu_1$. An application of the integration by parts formula from Remark \ref{IBP_formel_ad} results in
\begin{align*}
\int_U \alpha f^2+(CD f,Df)_U\dm_1= \int_U g f\dm_1.
\end{align*} 
To show the second equation we differentiate \eqref{equationUreg} with respect to the $k$-th direction yielding 
\begin{align}\label{partialg}
 \alpha \partial_k f-N_0\partial_kf+(d_k,Q_1^{-1}CD f)_U=\partial_k g
\end{align}Now we multiply the equation above with $\partial_l f(d_k,Cd_l)_U$. In order to structure the arguments we treat the resulting terms separately. If we sum over all indices's a direct calculation shows that the first and third term on the left hand side of the equation above is equal to $\alpha (CD f,D f)_U$ and $(CD f,Q_1^{-1}CD f)_U$, respectively. The right hand side of the equation is then equal to $(D g,CDf)_U$. We also get
\begin{align*}
\sum_{k,l=1}^{\infty}(d_k,Cd_l)_U\int_U -N_0\partial_k f \partial_l f\dm_1&=\sum_{k,l=1}^{\infty}(d_k,Cd_l)_U\int_U (CD \partial_k f, D \partial_l f)_U\dm_1\\
&=\int_U \mathrm{tr}[(CD^2f)^2]\dm_1,
\end{align*}and therefore the second equation from the statement is shown. Rearranging the terms of the equation we just derived yields
\begin{align*}
 \int_U\norm{Q_1^{-\frac{1}{2}}CDf}^2_U+\mathrm{tr}[(CD^2f)^2]\dm_1=\int_U (D(-N_0f),CDf)_U\dm_1.
\end{align*}Now it holds
\begin{multline*}
\int_U (D(-N_0f),CDf)_U\dm_1=\sum_{k,l=1}^{\infty}(d_k,Cd_l)_U\int_U \partial_k (-N_0f)\partial_lf\dm_1\\
=\sum_{k,l=1}^{\infty}(d_k,Cd_l)_U\int_U N_0f(\partial_{kl}f-(u,Q_1^{-1}d_k)_U\partial_lf)\dm_1=\int_U (N_0f)^2\dm_1,
\end{multline*}where we used that $N_0f\in W^{1,2}(U,\mu_1,\R)$ and Remark \ref{IBP_formel_ad}.

Note that the infinite sums in the calculations above are actually finite ones.
\end{proof}
\end{theorem}
\begin{remark}
For $f\in \fc B_U)$ there is some $n\in\N$, s.t. $Df\in B_U^n$. In particular $CDf\in B_U^n$, by Hypothesis \ref{hypoCadv}. Therefore
\begin{align}\label{coercive}
\frac{1}{\lambda_1}\norm{CDf}^2_U\leq \norm{Q_1^{-\frac{1}{2}}CDf}^2_U.
\end{align}Hence by the last equality in Theorem \ref{regestimates} we achieve
\begin{align}\label{eqreg}
 \begin{split}
 \int_U\frac{1}{\lambda_1}\norm{CDf}^2_U+\mathrm{tr}[(CD^2f)^2]\dm_1&\leq\int_U\norm{Q_1^{-\frac{1}{2}}CDf}^2_U+\mathrm{tr}[(CD^2f)^2]\dm_1\\
 &=\int_U (N_0f)^2\dm_1.
 \end{split}
\end{align}
\end{remark}
\begin{corollary}\label{domaindl}It holds $D(N_0)\subset W^{1,2}_{Q_1^{-\frac{1}{2}}C}(U,\mu_1,\R)\cap W^{2,2}_{C}(U,\mu_1,\R)$ and for all $f\in D(N_0)$ we have
\begin{align}\label{reg_DN0}
\begin{split}
 \int_U\frac{1}{\lambda_1}\norm{CDf}^2_U+\mathrm{tr}[(CD^2f)^2]\dm_1&\leq\int_U\norm{Q_1^{-\frac{1}{2}}CDf}^2_U+\mathrm{tr}[(CD^2f)^2]\dm_1\\
 &=\int_U (N_0f)^2\dm_1.
 \end{split}
\end{align}
\begin{proof}
Given $f\in D(N_0)$. We find a sequence $(f_n)_{n\in\N}\subset \fc B_U)$ such that $f_n\rightarrow f$ and $N_0f_n\rightarrow N_0f$ in $L^2(U,\mu_1,\R)$ as $n\rightarrow \infty$. By (In)equality \eqref{eqreg} $(f_n)_{n\in\N}$ is a Cauchy-sequence in $ W^{1,2}_{Q_1^{-\frac{1}{2}}C}(U,\mu_1,\R)$ as well as in $W^{2,2}_{C}(U,\mu_1,\R)$. Hence $f\in W^{1,2}_{Q_1^{-\frac{1}{2}}C}(U,\mu_1,\R)\cap W^{2,2}_{C}(U,\mu_1,\R)$ and the (in)equality of the statement is shown.
\end{proof}
\end{corollary}

Using (In)equality \eqref{reg_DN0} and Neumann's theorem we are able to deal with perturbations of $(N_0,\fc B_U))$ as described in the following theorem.
\begin{theorem}\label{essdispot}Assume that $\Phi$ is in $W^{1,2}(U,\mu_1,\R)$, bounded from below and with $\norm{D\Phi}^2_{L^{\infty}(\mu_1)}< \frac{1}{4\lambda_1}$.	
The operator $(N,\fc B_U))$ defined by
\begin{align*}
\fc B_U)\ni f\mapsto Nf=&\mathrm{tr}[CD^2f]-(u,Q^{-1}CD f)_U\\
&-(D\Phi,CDf)_U\in L^2(U,\mu_1^{\Phi},\R),
\end{align*}
\begin{enumerate}
\item fulfills 
\begin{align*}
(Nf,g)_{L^2(U,\mu_1^{\Phi},\R)}=\int_U -(CDf,Dg)_U\dm_1^{\Phi},
\end{align*} for all $f,g\in \fc B_U)$, in particular is dissipative in $L^2(U,\mu_1^{\Phi},\R)$. Furthermore we have
\item the dense range condition $\overline{(Id-N)(\fc B_U))}=L^2(U,\mu_1^{\Phi},\R)$.
\end{enumerate}In particular $(N,\fc B_U))$ is essentially m-dissipative in $L^2(U,\mu_1^{\Phi},\R)$. The resolvent in $\alpha\in (0,\infty)$ of the closure $(N,D(N))$ is denoted by $R(\alpha,N)$.
\begin{proof}
The first item of the statement follows by the integration by parts formula from Corollary \ref{IBP_potential} together with the invariance properties of the involved operators.
For $f\in L^2(U,\mu_1,\R)$ set
\begin{align*}
Tf=-(D\Phi,CDR(1,N_0)f)_U.
\end{align*}Using the Cauchy-Schwarz inequality, Inequality \eqref{reg_DN0} and the assumption on $\Phi$ we observe
\begin{align*}
\norm{Tf}_{L^2(\mu_1)}^2&=\int_U (D\Phi,CDR(1,N_0)f)_U^2\dm_1\\
&\leq \norm{D\Phi}^2_{L^{\infty}(\mu_1)}\int_U \norm{CDR(1,N_0)f}^2_U\dm_1<\frac{1}{4}\int_U (N_0R(1,N_0)f)^2\dm_1\\
&=\frac{1}{4}\int_U (f-R(1,N_0)f)^2\dm_1\leq \norm{f}_{L^2(\mu_1)}^2.
\end{align*}Therefore the linear operator $T:L^2(U,\mu_1,\R)\rightarrow L^2(U,\mu_1,\R)$ is well-defined with operator norm less than one. Hence by the Neumann-Series theorem we obtain that $(Id-T)^{-1}$ exists in $\mathcal{L}(L^2(U,\mu_1,\R))$. In particular for a given $g\in L^2(U,\mu_1,\R)$ we find $f\in L^2(U,\mu_1,\R)$ with $f-Tf=g$ in $L^2(U,\mu_1,\R)$. Since $(N_0,D(N_0))$ is m-dissipative, there is $h\in D(N_0)$ with $(Id-N_0)h=f$. This yields
\begin{align*}
(Id-N_0)h+(D\Phi,CDh)_U=f+(D\Phi,CDR(1,N_0)f)_U=f-Tf=g,
\end{align*}i.e. $L^2(U,\mu_1,\R)\subset(Id-N)(D(N_0))$.
Since $L^2(U,\mu_1,\R)$ is dense in $L^2(U,\mu_1^{\Phi},\R)$ it is left to show that the closure in $L^2(U,\mu_1^{\Phi},\R)$ of $(Id-N)(\fc B_U))$ contains $(Id-N)(D(N_0))$. This is true. Indeed for given $f\in D(N_0)$ we find a sequence $(f_n)_{n\in\N}\subset \fc B_U)$ s.t.~$f_n\rightarrow f$ and $N_0f_n\rightarrow N_0f$ in $L^2(U,\mu_1,\R)$. In view of the assumption on $D\Phi$, the Inequality \eqref{reg_DN0} and the fact that $\rho_{\Phi}=\frac{1}{c_{\Phi}}e^{-\Phi}$ is bounded we can estimate
\begin{align*}
&\norm{(Id-N)f-(Id-N)f_n}^2_{L^2(\mu_1^{\Phi})}=\norm{N(f-f_n)}^2_{L^2(\mu_1^{\Phi})}\\
&\leq2\norm{N_0(f-f_n)}^2_{L^2(\mu_1^{\Phi})}+2\int_U (D\Phi,CD(f-f_n))_U^2\dm_1^{\Phi}\\
&\leq 2\norm{\rho_{\Phi}}_{L^{\infty}(\mu_1)}\big(\norm{N_0(f-f_n)}^2_{L^2(\mu_1)}+\norm{D\Phi}^2_{L^{\infty}(\mu_1)}\int_U \norm{CD(f-f_n)}^2_U\dm_1\big)\\
& \leq 2\norm{\rho_{\Phi}}_{L^{\infty}(\mu_1)}\big(\norm{N_0(f-f_n)}^2_{L^2(\mu_1)}+\frac{1}{4}\norm{N_0(f-f_n)}^2_{L^2(\mu_1)}\big).
\end{align*}Using this estimate we conclude  
\begin{align*}
(Id-N)(D(N_0))\subset \overline{(Id-N)(\fc B_U))},
\end{align*}where the closure above is taken in $L^2(U,\mu_1^{\Phi},\R)$. Therefore the claim is shown.
\end{proof}
\end{theorem}

The following lines are devoted to derive a $L^2(U,\mu_1^{\Phi},\R)$ regularity estimate for the first and second order derivatives of a function $f\in\fc B_U)$ in terms of $g\in L^2(U,\mu_1,\R)$ related via
\begin{align}\label{equationUregpotential}
 \alpha f-Nf=g,
\end{align}for some given $\alpha \in (0,\infty)$.

\begin{theorem}\label{regestimatespotential}Assume that $\Phi:U\rightarrow \R$ is Fr\'{e}chet differentiable, bounded from below and $D\Phi:U\rightarrow U$ is Lipschitz continuous.
If $f\in \fc B_U)$, $g$ and $\alpha\in (0,\infty)$ are as in Equation \eqref{equationUregpotential} it holds
\begin{align*}
&\int_U\alpha f^2+(CD f,D f)_U\dm_1^{\Phi}=\int_Ugf\dm_1^{\Phi}\\
 &\int_U\alpha (CD f,D f)_U+\norm{Q_1^{-\frac{1}{2}}CDf}^2_U+\mathrm{tr}[(CD^2f)^2]+(D^2\Phi CDf,CDf)_U\dm_1^{\Phi}\\
 &=\int_U (D g,CD f)_U\dm_1^{\Phi}.
\end{align*}In particular we have
\begin{align*}
&\int_U	\norm{Q_1^{-\frac{1}{2}}CDf}^2_U+\mathrm{tr}[(CD^2f)^2]+(D^2\Phi CDf,CDf)_U\dm^{\Phi}
 =\int_U (Nf)^2\dm_1^{\Phi}.
\end{align*}
\begin{proof}The first equation follows by multiplying \eqref{equationUregpotential} with $f$, an integration over $U$ with respect to $\mu_1^{\Phi}$ and an application of the first item in Theorem \ref{essdispot}.
To show the second equation we differentiate \eqref{equationUregpotential} with respect to the $k$-th direction resulting in 
\begin{align*}
 \alpha \partial_k f-N\partial_kf+(d_k,Q_1^{-1}CD f)_U+\sum_{i=1}^{\infty} (d_i,CDf)_U\partial_{ki}\Phi=\partial_k g.
\end{align*}Note that the infinite sum in the line above is actually a finite one. Moreover $\partial_{ki}\Phi$ exists $\mu_1$-a.e., since the Lipschitz continuous function $\partial_i\Phi :U\rightarrow \R$ is Gateaux differentiable $\mu_1$-a.e.~by \cite[Proposition 10.11]{DaPrato_intro_infinite_dim_analysis}.
Now we multiply the equation above with $\partial_l f(d_k,Cd_l)_U$. If we sum over all indices's a direct calculation shows that the first and third term as well as the right hand side is equal to $\alpha (CD f,D f)_U$, $\norm{Q_1^{-\frac{1}{2}}CDf}^2_U$ and $(D g,CDf)_U$, respectively. For the second term we calculate
\begin{align*}
 &\sum_{k,l=1}^{\infty}(d_k,Cd_l)_U\int_U -N\partial_k f \partial_l f\dm_1^{\Phi}\\
 =&\sum_{k,l=1}^{\infty}(d_k,Cd_l)_U\int_U (CD \partial_k f, D \partial_l f)_U\dm_1^{\Phi}=\int_U \mathrm{tr}[(CD^2f)^2]\dm_1^{\Phi}.
\end{align*}Moreover we have
\begin{align*}
\sum_{k,l,i=1}^{\infty}(d_i,CDf)_U\partial_lf(d_k,Cd_l)_U\partial_{ki}\Phi &=\sum_{k,i=1}^{\infty} (d_i,CDf)_U(d_k,CDf)_U\partial_{ki}\Phi\\
&=(D^2\Phi CDf,CDf)_U,
\end{align*}from which we conclude the second equation. As in Theorem \ref{regestimates} we can rearrange the terms of the second equation to get
\begin{align*}
&\int_U\norm{Q_1^{-\frac{1}{2}}CDf}^2_U+\mathrm{tr}[(CD^2f)^2]+(D^2\Phi CDf,CDf)_U\dm_1^{\Phi}\\
 &=\int_U (D(-Nf),CDf)_U\dm_1^{\Phi}.
\end{align*}
Note that $N_0f,(D\Phi,CDf)_U\in W^{1,2}(U,\mu_1,\R)$, by \cite[Proposition 10.11]{DaPrato_intro_infinite_dim_analysis} and \cite[Proposition 10.9]{DaPrato_intro_infinite_dim_analysis} as well as $\partial_jf\rho_{\Phi}\in W^{1,2}(U,\mu_1,\R)$, $j\in\N$. Using the integration by parts formula from Remark \ref{IBP_formel_ad} we therefore get
\begin{align*}
&\int_U (D(-Nf),CDf)_U\dm_1^{\Phi}\\
&=\sum_{i,j=1}^{\infty}(Cd_i,d_j)_U\int_U \big(\partial_i(-N_0f)+ \partial_i(D\Phi,CDf)_U\big)\partial_jf\rho_{\Phi}\dm_1\\
&=\sum_{i,j=1}^{\infty}(Cd_i,d_j)_U\int_U Nf\big(\partial_{i}(\partial_jf\rho_{\Phi})-(u,Q_1^{-1}d_i)_U\partial_jf\rho_{\Phi}\big)\dm_1=\int_U (Nf)^2\dm_1^{\Phi}.
\end{align*}Note that the infinite sums in the calculations above are actually finite ones.
\end{proof}
\end{theorem}
\begin{remark}\label{reg_g_potential}Suppose we are in the situation of Theorem \ref{regestimatespotential}. Using the Cauchy-Schwarz inequality and the first equation in Theorem \ref{regestimatespotential} we obtain
\begin{align*}
\int_U \alpha f^2+(CDf,Df)_U\dm_1^{\Phi}&\leq\norm{g}_{L^2(\mu_1^{\Phi})}\norm{f}_{L^2(\mu_1^{\Phi})}=\norm{g}_{L^2(\mu_1^{\Phi})}\norm{R(\alpha,N)g}_{L^2(\mu_1^{\Phi})}\\
&\leq \frac{1}{\alpha} \norm{g}_{L^2(\mu_1^{\Phi})}^2.
\end{align*}
Now additionally suppose that $\Phi$ is a convex function. Hence we can estimate using the third equation in Theorem \ref{regestimatespotential} and the convexity of $\Phi$
\begin{multline}
\int_U \mathrm{tr}[(CD^2f)^2]+\norm{Q_1^{-\frac{1}{2}}CDf}_U^2\dm_1^{\Phi}\leq \int_U (Nf)^2\dm_1^{\Phi}= \int_U(\alpha f-g)^2\dm_1^{\Phi}\\
\leq 2\int_U(\alpha f)^2+g^2\dm_1^{\Phi}\leq 4\int_U g^2\dm_1^{\Phi}.
\end{multline}
\end{remark}
\begin{hypo}\label{hypo_potential_convex}$ $
\begin{enumerate}
\item The potential $\Phi$ is in $\in W^{1,2}(U,\mu_1,\R)$, convex, bounded from below and lower semicontinuous.
\item $\int_U \norm{D\Phi}_{U}^p\dm_1< \infty$ for some $p\in (2,\infty)$.
\end{enumerate}
\end{hypo}
\begin{remark}\label{yoshidaapp}For a potential $\Phi$ fulfilling Hypothesis \ref{hypo_potential_convex}  one can introduce the so called Yoshida approximation $\Phi_{t}$, $t>0$, defined by
\begin{align*}
\Phi_{t}(u)=\inf_{x\in U}\big\lbrace \Phi(x)+\frac{\norm{u-x}_U^2}{2t}\big\rbrace.
\end{align*}One can show that for all ${t>0}$ is the Yoshida approximation $\Phi_{t}:U\rightarrow (-\infty,\infty]$ is convex and Fr\'{e}chet differentiable with
\begin{enumerate}
\item $-\infty<\inf_{x\in U}\Phi(x)\leq\Phi_{t}(u)\leq \Phi(u)$ for all $u\in U$,
\item $\lim_{t\rightarrow 0}\Phi_{t}(u)=\Phi(u)$ for all $u\in U$,
\item $\norm{D\Phi_{t}(u)}_U\leq \norm{D\Phi(u)}_U$ for $\mu_1$-a.e. $u\in U$ and
\item $\lim_{t\rightarrow 0}D\Phi_{t}(u)=D\Phi(u)$ for $\mu_1$-a.e. $u\in U$.
\end{enumerate}Furthermore $D\Phi_{t}$ is Lipschitz continuous for all $t >0$. A proof of these statements can be found in \cite{yoshida}.
\end{remark}
\begin{theorem}\label{theorem_reg}
Suppose $\Phi$ fulfills Hypothesis \ref{hypo_potential_convex}. Then for $f\in \fc B_U)$ and $g=\alpha f-Nf$, $\alpha\in(0,\infty)$, as in Equation \eqref{equationUregpotential} we have
\begin{align*}
\int_U \alpha f^2+(CDf,Df)_U\dm_1^{\Phi}&\leq\frac{1}{\alpha} \int_U g^2\dm_1^{\Phi},\\
\int_U \mathrm{tr}[(CD^2f)^2]+\norm{Q_1^{-\frac{1}{2}}CDf}_U^2\dm_1^{\Phi}&\leq 4\int_U g^2\dm_1^{\Phi}.
\end{align*}
\begin{proof}
Let $(\Phi_{t})_{t>0}$ be the Yoshida approximation of $\Phi$. For $t>0$, define $g_{t}$ by
\begin{align*}
     \alpha f-\mathrm{tr}[CD^2f]+(\cdot,Q_1^{-1}CD f)_U-(D\Phi_{t},CDf)_U=g_{t}.
\end{align*}By Remark \ref{reg_g_potential} we obtain
\begin{align*}
\int_U (\alpha f^2+(CDf,Df)_U)\rho_{\Phi_{t}}\dm_1&\leq \frac{1}{\alpha}\int_U g_{t}^2\rho_{\Phi_{t}}\dm_1,\\
\int_U (\mathrm{tr}[(CD^2f)^2]+\norm{Q_1^{-\frac{1}{2}}CDf}_U^2) \rho_{\Phi_{t}}\dm_1&\leq 4\int_U g_{t}^2 \rho_{\Phi_{t}}\dm_1.
\end{align*}
By Remark \ref{yoshidaapp}, $\rho_{\Phi_{t}}=\frac{1}{c_{\Phi_{t}}}e^{-\Phi_{t}}$ is bounded by a constant $\theta$ independent of $t$. In particular $(\mathrm{tr}[(CD^2f)^2]+\norm{Q_1^{-\frac{1}{2}}CDf}_U^2)\rho_{\Phi_{t}}$ and $(\alpha f^2+(CDf,Df)_U)\rho_{\Phi_{t}}$ are bounded by $(\mathrm{tr}[(CD^2f)^2]+\norm{Q_1^{-\frac{1}{2}}CDf}_U^2)\theta$ and $(\alpha f^2+(CDf,Df)_U)\theta$, respectively. Since they converge pointwisely to $(\mathrm{tr}[(CD^2f)^2]+\norm{Q_1^{-\frac{1}{2}}CDf}_U^2) \rho_{\Phi}$ and $(\alpha f^2+(CDf,Df)_U)\rho_{\Phi}$ we know that the left hand sides of the inequalities above converge to $\int_U (\mathrm{tr}[(CD^2f)^2]+\norm{Q_1^{-\frac{1}{2}}CDf}_U^2)\dm_1^{\Phi}$ and $\int_U \alpha f^2+(CDf,Df)_U\dm_1^{\Phi}$, respectively. It also holds 
\begin{align*}
\abs{\int_U g_{t}^2 \rho_{\Phi_{t}}-g^2\rho_{\Phi}\dm_1}&\leq\abs{\int_U g_{t}^2 (\rho_{\Phi_{t}}-\rho_{\Phi})\dm_1}+\abs{\int_U (g_{t}^2-g^2)\rho_{\Phi}\dm_1}\\
&=\abs{\int_U g_{t}^2 (e^{-\Phi_{t}}-e^{-\Phi})\dm_1}+\abs{\norm{g_{t}}^2_{L^2(\mu_1^{\Phi})}-\norm{g}^2_{L^2(\mu_1^{\Phi})}}.
\end{align*}Note that $g_{t}^2$ can be bounded independent of $t$ by an $\mu_1$-integrable function, hence the first term in the above inequality goes to zero as $t$ goes to zero by another application of the dominated convergence theorem. The second term also tends to zero. Indeed the Cauchy-Schwarz inequality and the definitions of $g$ and $g_{t}$ yields
\begin{align*}
\int_U (g-g_{t})^2\dm_1^{\Phi}&=\int_U(D\Phi-D\Phi_{t},CDf)_U^2 \dm_1^{\Phi}\\
&\leq \int_U \norm{D\Phi-D\Phi_{t}}_U^2\norm{CDf}_U^2\dm_1^{\Phi}.
\end{align*}Invoking the third and the fourth item of Remark \ref{yoshidaapp} and another application of the dominated convergence theorem yields that $g_{t}$ converges to $g$ in $L^2(U,\mu_1^{\Phi},\R)$ as $t$ goes to zero. In particular the corresponding norms in $L^2(U,\mu_1^{\Phi},\R)$ converge. All together this finishes the proof.
\end{proof}
\end{theorem}

Note that the regularity estimate we derived in Theorem \ref{theorem_reg} relies on Hypothesis \ref{hypo_potential_convex} which is less restrictive than the assumptions in Theorem \ref{essdispot}. 
   
\section{The infinite-dimensional Langevin operator}\label{section_lan}
The essential m-dissipativity of finite-dimensional Langevin operators have been extensively studied in \cite{GrothausConrad1} and \cite{DissConrad} for singular potentials and even in a manifold setting in \cite{LangevinManifold}. In this section we want to extend these result to an infinite-dimensional setting, where as in the above references the non-sectorality of $L_{\Phi}$ causes difficulties.

As in the introduction, we fix two real separable Hilbert spaces $(U,(\cdot,\cdot)_U)$ and $(V,(\cdot,\cdot)_V)$ and consider the real separable Hilbert space $(W,(\cdot,\cdot)_W)$ defined by $W=U\times V$ and
\begin{align*}
((u_1,v_1),(u_2,v_2))_W=(u_1,u_2)_U+(v_1,v_2)_V,\quad (u_1,v_1),(u_2,v_2)\in W.
\end{align*}By $\mu_1$ and $\mu_2$ we denote two centered non-degenerate Gaussian measures on $(U,\mathcal{B}(U))$ and $(V,\mathcal{B}(V))$, respectively. The corresponding covariance operators are denoted by $Q_1\in \mathcal{L}_1^+(U)$ and $Q_2\in \mathcal{L}_1^+(V)$. We also fix two orthonormal basis $B_U=\di$ and $B_V=\ei$ of eigenvectors with corresponding eigenvalues $(\lambda_i)_{i\in\N}$ and $(\nu_i)_{i\in\N}$ of $Q_1$ and $Q_2$, respectively. W.l.o.g. we assume that $(\lambda_i)_{i\in\N}$ and $(\nu_i)_{i\in\N}$ are decreasing to zero. Furthermore  we set $B_W=(B_U,B_V)$.\\
As in Definition \ref{spaces} one can consider the orthogonal projections to $B_U^n$ and $B_V^n$, $n\in\N$. To avoid an overload of notation we omit to indicate if we project to $B_U^n$ and $B_V^n$ as it is clear from the context.

On $(W,\mathcal{B}(W))$ we consider the product measure $\mu=\mu_1\otimes\mu_2$. Using the separability of $U$ and $V$, \cite[Lemma 1.2]{FoundationsofModernProbability} it holds $\mathcal{B}(W)=\mathcal{B}(U)\otimes \mathcal{B}(V)$. Applying \cite[Theorem 1.12]{DaPrato_intro_infinite_dim_analysis} one can check that $\mu$ is a non-degenerate centered Gaussian measure with centered non-degenerate covariance operator $Q\in \mathcal{L}_1^+(W)$ defined by
\begin{align*}
W\ni (u,v)\mapsto Q(u,v)=(Q_1u,Q_2v)\in W.
\end{align*}\begin{definition}\label{smooth}
In $L^2(\mu)$ we denote by $\fc B_W)$ the space of finitely based smooth and bounded functions on $W$ defined by
\begin{align*}
&\fc B_W)\\
&=\lbrace W\ni (u,v)\mapsto \varphi(P_m(u),P_m(v))\in \R\mid m\in\N, \; \varphi\in C_b^{\infty}(\R^{m}\times \R^m)\rbrace
\end{align*} and correspondingly the space of finitely based smooth and bounded functions on $W$ only dependent on the first $n$ directions by
\begin{align*}
\fc B_W,n)=\lbrace W\ni (u,v)\mapsto \varphi(P_n(u),P_n(v))\in \R\mid \; \varphi\in C_b^{\infty}(\R^{n}\times \R^n)\rbrace.
\end{align*}
\end{definition}
Concerning derivatives of sufficient smooth functions $f:W\rightarrow \R$ recall the explanation in Remark \ref{derivatives}.
We set $D_1f=\sum_{i=1}^{\infty}(Df,(d_i,0))_Wd_i\in U$ and $D_2=\sum_{i=1}^{\infty}(Df,(0,e_i))_We_i\in V$ as well as $\partial_{i,1}f=(D_1f,d_i)_U$ and $\partial_{i,2}f=(D_2f,e_i)_V$. In particular we have
\begin{align*}
Df=\sum_{i=1}^{\infty}(Df,(d_i,0))_W(d_i,0)+\sum_{i=1}^{\infty}(Df,(0,e_i))_W(0,e_i)=(D_1f,D_2f).
\end{align*}Analogously we define $D_1^2f$, $D_2^2f$ as well as $\partial_{ij,1}f$ and $\partial_{ij,2}f$.

For given $n\in\N$, recall the image measures $\mu_1^n$ and $\mu_2^n$ from Lemma \ref{imagemeasure} w.r.t. $B_U$ and $B_V$, respectively and set $\mu^n=\mu_1^n\otimes \mu_2^n$ on $(\R^{n}\times \R^n,\mathcal{B}(\R^n)\otimes \mathcal{B}(\R^n))$. We also consider the Hilbert space
\begin{align*}
L^2_{B_{W_n}}(\mu)=\lbrace W\ni (u,v)\mapsto f(P_n(u),P_n(v))\in\R\mid f\in L^2(\mu^n)\rbrace,
\end{align*}similarly to Definition \ref{spaces}.
\begin{remark}\label{densitytwo}Arguing as in Lemma \ref{density} one can show that $\fc B_W)$ and $\fc B_W,n)$ are dense in $L^2(\mu)$ and $L^2_{B_{W_n}}(\mu)$, respectively.
\end{remark}
Moreover we fix operators $K_{12}\in \mathcal{L}(U;V)$, $K_{21}\in \mathcal{L}(V;U)$ and $K_{22}\in \mathcal{L}^+(U)$. Last but not least we consider a measurable potential $\Phi:U\rightarrow (-\infty,\infty]$ which is bounded from below and recall the measures $\mu_1^{\Phi}$ and $\mu^{\Phi}$. During the whole section we will assume the following hypothesis.
\begin{hypo}\label{stass}$ $
\begin{enumerate}
\item $K_{22}$ is positive and for all $n\in \N$ the space $B_V^n$ is  $K_{22}$ invariant. 
\item  $K_{12}^*=K_{21}$.
\item For all $n\in N$ we have $K_{12}(B_U^n)\subset B_V^n$ and  $K_{21}(B_V^n)\subset B_U^n$.
\item The operator $(-K_{22}Q_2^{-1},D(K_{22}Q_2^{-1}))$ on $V$ generates a strongly continuous contraction semigroup $(\exp(-tK_{22}Q_2^{-1}))_{t\geq 0}$ and 
there is a $\tau_2 > 0$ such that 
\begin{align*}
(-K_{22}Q_2^{-1}v,v)_V\leq -\tau_2\norm{v}_V^2,\quad v\in D(K_{22}Q_2^{-1})= D(Q_2^{-1}).
\end{align*}
\item $\Phi\in W^{1,2}(U,\mu_1,\R)$.
\end{enumerate}
\end{hypo}
Next we define the infinite-dimensional Langevin operator $(L_{\Phi},\fc B_W))$. We will realize in Remark \ref{decompose}, that the invariance properties of $K_{12}$, $K_{21}$ and $K_{22}$ included in the hypothesis above, ensures that $(L_{\Phi},\fc B_W))$ has a useful decomposability property.
\begin{definition}\label{rigLangevinop}

We define $(L_{\Phi},\fc B_W))$ in $L^2(\mu^{\Phi})$ by 
\begin{align*}
\fc B_W)  \ni f \mapsto L_{\Phi}f=S_{\Phi}f-A_{\Phi}f\in L^2(\mu^{\Phi}),
\end{align*} where for $f\in \fc B_W)$, $S_{\Phi}f$ and $A_{\Phi}f$ are given by
\begin{align*}
S_{\Phi}f=&\mathrm{tr}[K_{22}D^2_2f]-(v,Q_2^{-1}K_{22}D_2f)_V ,\\
A_{\Phi}f=&(u,Q_1^{-1}K_{21}D_2f)_U+(D \Phi(u),K_{21}D_2f)_U-(v,Q_2^{-1}K_{12}D_1f)_V.
\end{align*}
The designation of $S_{\Phi}$ and $A_{\Phi}$ is not accidental, as we see show in the next lemma, that $(S_{\Phi},\fc B_W))$ is symmetric and $(A_{\Phi},\fc B_W))$ antisymmetric.

Remember that $\fc B_W)$ is dense in $L^2(\mu_1^{\Phi})$ by Remark \ref{densitytwo} and expressions as $Q_2^{-1}K_{22}D_2f$, $Q_2^{-1}K_{22}D_2f$ and  $Q_1^{-1}K_{21}D_2f$ are reasonable due to Hypothesis \ref{stass} and Remark \ref{domain}.
\end{definition}

Using the integration by parts formula from Theorem  \ref{IBP_formel} together with the invariance properties of $K_{22}$, $K_{21}$ and $K_{12}$ one can derive the following important lemma.

\begin{lemma}\label{propL}It holds
\begin{enumerate}
\item[(i)]$(S_{\Phi},\fc B_W))$ is symmetric and dissipative in $L^2(\mu^{\Phi})$.
\item[(ii)]$(A_{\Phi},\fc B_W))$ is antisymmetric in $L^2(\mu^{\Phi})$.
\item[(iii)]$1\in \fc B_W)$ with $L_{\Phi}1=0$ and in particular $\mu^{\Phi}$ is invariant for $(L_{\Phi},\fc B_W))$ in the sense that
\begin{align*}
\int_W L_{\Phi}f\mathrm{d}\mu^{\Phi}=0\quad\text{for\;all}\quad f\in  \fc B_W).
\end{align*}
\item[(iv)]$(L_{\Phi},\fc B_W))$ is dissipative in $L^2(\mu^{\Phi})$ and for all $f,g\in \fc B_W)$ it holds
\begin{align*}
-&\int_W L_{\Phi}f g\mathrm{d}\mu^{\Phi}\\
=&\int_W (D_2f,K_{22}D_2g)_V-(D_1f,K_{21}D_2g)_U+(D_2f,K_{12}D_1g)_V \mathrm{d}\mu^{\Phi}.
\end{align*}
\end{enumerate}
\end{lemma}
By \cite[Proposition 3.14]{OSG0} densely defined dissipative operators are closable. Since $(L,\fc B_W))$, $(S_{\Phi},\fc B_W))$ and $(A,\fc B_W))$ are densely defined dissipative operators in $L^2(\mu^{\Phi})$, it is reasonable to denote their closures by $(L_{\Phi},D(L_{\Phi}))$, $(S_{\Phi},D(S_{\Phi}))$ and $(A_{\Phi},D(A_{\Phi}))$.
The overall goal is to show essential m-dissipativity of $(L_{\Phi},\fc B_W))$, i.e.~m-dissipativity of $(L_{\Phi},D(L_{\Phi}))$ in $L^2(\mu^{\Phi})$. This will be done using the strategy described below.
\begin{enumerate}
\item At first we consider $\Phi=0$. In this case we set $L=L_0$, $S=S_0$ and $A=A_0$. We decompose $L$ into countable many finite-dimensional Langevin operators and translate our problem into a finite-dimensional one. A unitary transformation is used to get nice representations $S_n$ $A_n$ and $L_n$  of $S$,$A$ and $L$ in $L^2(\R^n\times\R^n,\mu^n,\R)$. In Proposition \ref{prop_ess_S_N} we show essential m-dissipativity of $S_n$, which is defined below. 
\item We follow the strategy of \cite{GrothausConrad1}, to obtain essential m-dissipativity of $L_n$. 
\item Instead of assuming that $\Phi=0$ we consider a measurable potential $\Phi :U\mapsto (-\infty,\infty]$, which is bounded from below fulfilling the assumptions as asserted in Theorem \ref{ess_diss_Lphi}. Invoking a similar strategy as in Theorem \ref{essdispot} we get the essential m-dissipativity of $(L_{\Phi},\fc B_W))$ in $L^2(\mu^{\Phi})$.
\end{enumerate}
As described above we start with $\Phi=0$, hence the infinite-dimensional measure $\mu^{\Phi}$ reduces to the infinite-dimensional centered non-degenerate Gaussian measure $\mu$ with covariance operator $Q$.
\begin{remark}\label{decompose}
Given $n\in\N$. Set
\begin{align*}
K_{22,n}=((K_{22}e_i,e_j))_{ij=1}^n,\quad K_{12,n}=((K_{12}d_i,e_j))_{ij=1}^n,\quad K_{21,n}=K_{12,n}^*.
\end{align*}
We have for $f=\varphi(P_n(\cdot),P_n(\cdot))\in \fc B_W)$ and $(u,v)\in W$
\begin{align*}
Sf(u,v)=&\mathrm{tr}[K_{22,n}D^2_2\varphi(P_n(u),P_n(v))]-\langle P_n(v),Q_{2,n}^{-1}K_{22,n}D_2\varphi(P_n(u),P_n(v))\rangle,\\
Af(u,v)=&\langle P_n(u),Q_{1,n}^{-1}K_{21,n}D_2\varphi(P_n(u),P_n(v))\rangle\\
&-\langle P_n(v),Q_{2,n}^{-1}K_{12,n}D_1\varphi(P_n(u),P_n(v))\rangle.
\end{align*}
Hence, it is reasonable to consider the symmetric operator $(S_n,C_b^{\infty}(\R^{n}\times\R^n))$ and antisymmetric operator $(A_n,C_b^{\infty}(\R^{n}\times\R^n))$ in $L^2(\mu^n)$, defined for $\varphi\in C_b^{\infty}(\R^{n}\times\R^n)$ by
\begin{align*}
S_n\varphi=&\mathrm{tr}[K_{22,n}D^2_2\varphi]-\langle y,Q_{2,n}^{-1}K_{22,n}D_2\varphi\rangle,\\
A_n\varphi=&\langle x,Q_{1,n}^{-1}K_{21,n}D_2\varphi\rangle-\langle y,Q_{2,n}^{-1}K_{12,n}D_1\varphi\rangle.
\end{align*}On $C_b^{\infty}(\R^{n}\times\R^n)$ we set $L_n=S_n-A_n$ in analogy to the definition of $L$.
\end{remark}

Before we start with the first proposition we need the following notation.
If $E=E_1\times E_2$, where $E_1$ and $E_2$ are sets, and $f_1 :E_1 \rightarrow \R$, $f_1 :E_1 \rightarrow \R$ are functions, we denote by $f_1 \otimes f_2$ the function $E \ni (x,y) \rightarrow f_1(x)f_2(y)$. If $G_1$ and $G_2$ are linear spaces of functions on $E_1$ and $E_2$, respectively, we denote by $G_1 \otimes G_2$ the linear span of the set of all functions of the form $f_1\otimes f_2$, $f_2 \in G_1$, $f_2 \in G_2$.

\begin{proposition}\label{prop_ess_S_N}The operator $(S_n,C_b^{\infty}(\R^{n}\times\R^n))$ is essentially m-dissipative in $L^2(\mu^n)$. In particular its closure exists and can be denoted by $(S_n,D(S_n))$.
\begin{proof}
For arbitrary $\varphi\in C_b^{\infty}(\R^{n}\times\R^n)$ set $f=\varphi(P_n(\cdot),P_n(\cdot))\in \fc B_W)$. It holds
\begin{align*}
(S_n\varphi,\varphi)_{L^2(\mu^n)}=(Sf,f)_{L^2(\mu)}.
\end{align*}Hence the dissipativity of $(S_n,C_b^{\infty}(\R^n\times\R^n))$ follows by the dissipativity of $(S,\fc B_W))$. By the Lumer-Phillips theorem it is left to show that $(Id-S_n)(C_b^{\infty}(\R^{n}\times\R^n))$ is dense in $L^2(\mu^n)$. Given $\varepsilon>0$ and a function $\phi=\phi_1\otimes\phi_2\in  C_b^{\infty}(\R^n)\otimes C_b^{\infty}(\R^n)$. Assume w.l.o.g that $\norm{\phi_1}_{L^2(\mu_1^n)}>0$. Applying Theorem \ref{ornstein_ess} we know that there is a function $\varphi_2\in C_b^{\infty}(\R^n)$ such that
\begin{align*}
\norm{\varphi_2-\mathrm{tr}[K_{22,n}D^2\varphi_2]+\langle y,Q_{2,n}^{-1}K_{22,n}D\varphi_2\rangle-\phi_2}_{L^2(\mu_2^n)}< \frac{\varepsilon}{\norm{\phi_1}_{L^2(\mu_1^n)}}.
\end{align*}Set $\varphi=\phi_1\otimes\varphi_2\in  C_b^{\infty}(\R^{n}\times\R^n)$ to obtain
\begin{align*}
&\norm{(Id-S_n)\varphi-\phi}_{L^2(\mu^n)}\\
=&\norm{\phi_1}_{L^2(\mu_1^n)}\norm{\varphi_2-\mathrm{tr}[K_{22,n}D^2\varphi_2]+\langle y,Q_{2,n}^{-1}K_{22,n}D\varphi_2\rangle-\phi_2}_{L^2(\mu_2^n)}< \varepsilon.
\end{align*}Hence $C_b^{\infty}(\R^n)\otimes C_b^{\infty}(\R^n)$ is contained in the closure of $(Id-S_n)(C_b^{\infty}(\R^{n}\times\R^n))$. This is enough to show the claim as  $C_b^{\infty}(\R^n)\otimes C_b^{\infty}(\R^n)$ is dense in $L^2(\mu^n)$.
\end{proof}
\end{proposition}
Denote by $L^2_0(\lambda^n)$ the set of functions in $L^2(\lambda^{n})$ vanishing almost everywhere outside a bounded set in $\R^n$. By $\mathcal{S}(\R^n)$ we denote the set of rapidly decreasing functions from $\R^n$ to $\R$, also known as Schwartz functions.

\begin{definition}Define the unitary transformation $U:L^2(\lambda^n\otimes\lambda^n)\rightarrow L^2(\mu^n)$ by
\begin{align*}
f \mapsto f\frac{1}{c_1^n}\exp(\frac{1}{4}\langle Q^{-1}_{1,n}x,x\rangle) \frac{1}{c_2^n}\exp(\frac{1}{4}\langle Q^{-1}_{2,n}y,y\rangle)f
\end{align*}where $c_1^n$ ans $c_2^n$ are the corresponding normalizing constants.
\end{definition}

\begin{proposition}The operator $(\tilde{S}_n,L_0^2(\R^n)\otimes C_0^{\infty}(\R^n))$ in $L^2(\lambda^n\otimes\lambda^n)$, defined by $\tilde{S}_n=U^{-1}S_nU$ is essentially m-dissipative in $L^2(\lambda^n\otimes\lambda^n)$. For $\varphi\in L_0^2(\R^n)\otimes C_0^{\infty}(\R^n)$ it holds
\begin{align*}
\tilde{S}_n\varphi=&\mathrm{tr}[K_{22,n}D_2^2\varphi]-\frac{1}{4}\langle Q_{2,n}^{-1}K_{22,n}Q_{2,n}^{-1}y,y\rangle\varphi+\frac{\mathrm{tr}[K_{22,n}Q_{2,n}^{-1}]}{2}\varphi.
\end{align*}
\begin{proof}
Using a standard approximation argument one can show that $(C_0^{\infty}(\R^n)\otimes C_0^{\infty}(\R^n))$ is a core for $(S_n,D(S_n))$. As unitary transformations preserves essential m-dissipativity, we know that $(\tilde{S}_n,U^{-1}(C_0^{\infty}(\R^n)\otimes C_0^{\infty}(\R^n))$ is essential m-dissipative in $L^2(\lambda^n\otimes\lambda^n)$. Since $U^{-1}(C_0^{\infty}(\R^n)\otimes C_0^{\infty}(\R^n))=C_0^{\infty}(\R^n)\otimes C_0^{\infty}(\R^n)$ the same holds true for $(\tilde{S}_n,C_0^{\infty}(\R^n)\otimes C_0^{\infty}(\R^n))$. By construction $(\tilde{S}_n,L^2_0(\lambda^n)\otimes C_0^{\infty}(\R^n))$ is a dissipative extension of $(\tilde{S}_n,C_0^{\infty}(\R^n)\otimes C_0^{\infty}(\R^n))$ and therefore the first assertion is shown. The second assertion follows by a direct calculation. 
\end{proof}
\end{proposition}

Until otherwise stated all our functions in the upcoming considerations are assumed to be \textbf{complex valued}. Observe that the involved complex valued function spaces are isomorphic to their real counterparts. If an operator leaves the real valued functions invariant, essential m-dissipativity of this operator in the complex valued case implies the essential m-dissipativity in the real valued case, compare \cite[Lemma 1.1.21.]{DissConrad}.

In the next lemma we perturbate the operator $(\tilde{S}_n,L_0^2(\R^n)\otimes C_0^{\infty}(\R^n))$ with the operator $(B_n,L_0^2(\R^n)\otimes C_0^{\infty}(\R^n))$ given by
\begin{align*}
L_0^2(\R^n)\otimes C_0^{\infty}(\R^n)\ni \varphi\mapsto B_n\varphi=i\langle Q_{2,n}^{-1}y,K_{12,n}x\rangle\varphi \in L^2(\lambda^n\otimes\lambda^n),
\end{align*}without loosing essential m-dissipativity in $ L^2(\lambda^n\otimes\lambda^n)$. The key in the proof is \cite[Lemma 1.4.5.]{DissConrad}, invoking the concept of relative boundedness of operators and complete families of orthogonal projections compare \cite[Chapter 1.4]{DissConrad}. In our particular case we have a look at the complete family of orthogonal projections $(M_k)_{k\in\N}$ on $L^2(\lambda^n\otimes\lambda^n)$ given by the multiplication operators
\begin{align*}
L^2(\lambda^n\otimes\lambda^n) \ni \varphi\mapsto M_k\varphi=g_k\varphi\in L^2(\lambda^n\otimes\lambda^n),
\end{align*}where $g_k(x,y)=\mathbbm{1}_{[k−1,k]}(\abs{x})$ for all $(x,y)\in \R^n\times \R^n$.

\begin{lemma}\label{Diss_fourier_pert}$(\tilde{S}_n+B_n,L_0^2(\R^n)\otimes C_0^{\infty}(\R^n))$ is essential m-dissipative in $L^2(\lambda^n\otimes\lambda^n)$.
\begin{proof}
First of all note that $\langle Q_{2,n}^{-1}y,K_{12,n}x\rangle$ is real for all $(x,y)\in\R^n\times\R^n$, hence for the real part of $(B_n\varphi,\varphi)_{L^2(\lambda^n\otimes\lambda^n)}$, $\varphi\in L_0^2(\R^n)\otimes C_0^{\infty}(\R^n)$, it holds
\begin{align*}
\Re(B_n\varphi,\varphi)_{L^2(\lambda^n\otimes\lambda^n)}=0.
\end{align*}In particular $(B_{n},L_0^2(\R^n)\otimes C_0^{\infty}(\R^n))$ is dissipative in ${L^2(\lambda^n\otimes\lambda^n)}$.

By construction $M_k(L_0^2(\R^n)\otimes C_0^{\infty}(\R^n))\subset L_0^2(\R^n)\otimes C_0^{\infty}(\R^n)$, $M_kB_n=B_nM_k$ and $M_k\tilde{S}_n=\tilde{S}_n M_k$, for all $k\in \N$. Set $B_n^k=M_kB_n$ and $A_n^k=M_k\tilde{S}_n$ both with domain $M_k(L_0^2(\R^n)\otimes C_0^{\infty}(\R^n))$. In view of \cite[Lemma 1.4.5.]{DissConrad} it is left to show that each $B_n^k$ is $A_n^k$-bounded with $A_n^k$-bound less or equal to 1. Given $k\in\N$ and $\varphi\in M_k(L_0^2(\R^n)\otimes C_0^{\infty}(\R^n))$. Setting $\varphi_n=\abs{K_{12,n}x} \varphi$, $a=\abs{K_{22,n}^{-\frac{1}{2}}}$, $b=\abs{K_{12,n}}$ and $c=\mathrm{tr}[K_{22,n}Q_{2,n}^{-1}]$ we can estimate
\begin{align*}
&\norm{B_n\varphi}_{L^2(\lambda^n\otimes\lambda^n)}^2=\int_{\R^{n}\times\R^n}\langle Q_{2,n}^{-1}y,K_{12,n}x\rangle^2\varphi\overline{\varphi}\mathrm{d}\lambda^n\otimes\lambda^n\\
\leq& \int_{\R^{n}\times\R^n}\abs{K_{22,n}^{-\frac{1}{2}}}^2 \abs{K_{22,n}^{\frac{1}{2}}Q_{2,n}^{-1}y}^2\abs{K_{12,n}x}^2 \varphi\overline{\varphi}\mathrm{d}\lambda^n\otimes\lambda^n\\
=&a^2 \big(\abs{K_{22,n}^{\frac{1}{2}}Q_{2,n}^{-1}y}^2\varphi_n,\varphi_n\big)_{L^2(\lambda^n\otimes\lambda^n)}\\
\leq& a^2\big( -4 \mathrm{tr}[K_{22,n}D_2^2\varphi_n]+\abs{K_{22,n}^{-\frac{1}{2}}Q_{2,n}^{-1}y}^2\varphi_n-2c\varphi_n+2c\varphi_n,  \varphi_n\big)_{L^2(\lambda^n\otimes\lambda^n)}\\
=&4 a^2\big( (-\tilde{S}_n\varphi_n,\varphi_n\big)_{L^2(\lambda^n\otimes\lambda^n)}+\big(\frac{c}{2}\varphi_n, \varphi_n)_{L^2(\lambda^n\otimes\lambda^n)}\big)\\
\leq&  4 a^2 \norm{\abs{K_{12,n}x}^2 \varphi}_{L^2(\lambda^n\otimes\lambda^n)}\big( \norm{\tilde{S}_n\varphi}_{L^2(\lambda^n\otimes\lambda^n)}+\frac{c}{2}\norm{\varphi}_{L^2(\lambda^n\otimes\lambda^n)}\big).
\end{align*}Note that the estimation above relies on the Cauchy-Schwarz inequality, as well as the fact that the operator $\mathrm{tr}[K_{22,n}D_2^2]$ is dissipative in $L^2(\lambda^n\otimes\lambda^n)$.\\
Using that $\varphi$ is in $M_k(L_0^2(\R^n)\otimes C_0^{\infty}(\R^n))$ it holds
\begin{align*}
\norm{\abs{K_{12,n}x}^2 \varphi}_{L^2(\lambda^n\otimes\lambda^n)}\leq b^2k^2\norm{\varphi}_{L^2(\lambda^n\otimes\lambda^n)}.
\end{align*}The combination of these two estimates yields
\begin{align*}
\norm{B_n\varphi}_{L^2(\lambda^n\otimes\lambda^n)}^2\leq 4a^2   b^2k^2\big(\norm{\tilde{S}_{n}\varphi}_{L^2(\lambda^n\otimes\lambda^n)}\norm{\varphi}_{L^2(\lambda^n\otimes\lambda^n)}+\frac{c}{2}\norm{\varphi}_{L^2(\lambda^n\otimes\lambda^n)}^2\big).
\end{align*}Now the claim follows with \cite[Lemma 1.4.5.]{DissConrad}.
\end{proof}
\end{lemma}
\begin{remark}Note that $C_0^{\infty}(\R^n)\otimes C_0^{\infty}(\R^n)$ is a core for $(\tilde{S}_n+B_n,L_0^2(\R^n)\otimes C_0^{\infty}(\R^n))$. In particular $(\tilde{S}_n+B_n,C_0^{\infty}(\R^n)\otimes C_0^{\infty}(\R^n))$ and its dissipative extension $(\tilde{S}_n+B_n,\mathcal{S}(\R^n)\otimes C_0^{\infty}(\R^n))$ is essential m-dissipative in $L^2(\lambda^n\otimes\lambda^n)$.
Denote by $\mathcal{F}:L^2(\R^n,\lambda^n,\C)\rightarrow L^2(\R^n,\lambda^n,\C)$ the Fourier transformation. It is well known that $\mathcal{F}$ is unitary and it holds for all $\varphi\in S(\R^n)$ and multi-indices $s\in \N^n$
\begin{align}\label{Fourier}
\mathcal{F}(\partial_s\varphi)=i^{\abs{s}}x^{s}\mathcal{F}(\varphi).
\end{align}For $\varphi=\varphi_1\otimes \varphi_2\in \mathcal{S}(\R^n)\otimes C_0^{\infty}(\R^n)$, we set $\mathcal{F}_1\varphi=\mathcal{F}\varphi_1\otimes \varphi_2$ and extend $\mathcal{F}$ linearly to $\mathcal{S}(\R^n)\otimes C_0^{\infty}(\R^n)$. As for the classical Fourier transform we extend $\mathcal{F}_1:\mathcal{S}(\R^n)\otimes C_0^{\infty}(\R^n) \rightarrow L^2(\lambda^n\otimes\lambda^n)$ to a unitary operator operator from $L^2(\lambda^n\otimes\lambda^n)$ to $L^2(\lambda^n\otimes\lambda^n)$. Note that $\mathcal{F}(\mathcal{S}(\R^n))=\mathcal{S}(\R^n)$, hence also $\mathcal{F}_1(\mathcal{S}(\R^n)\otimes C_0^{\infty}(\R^n))=\mathcal{S}(\R^n)\otimes C_0^{\infty}(\R^n)$.\\
Now define  $\hat{S}_n=\mathcal{F}_1(\tilde{S}_n+B_n)\mathcal{F}_1^{-1}$ on $\mathcal{S}(\R^n)\otimes C_0^{\infty}(\R^n)$. An application of identity \eqref{Fourier} yields for all $\varphi\in \mathcal{S}(\R^n)\otimes C_0^{\infty}(\R^n)$
\begin{align*}
\hat{S}_n\varphi=&\mathrm{tr}[K_{22,n}D_2^2\varphi]-\frac{1}{4}\langle Q_{2,n}^{-1}K_{22,n}Q_{2,n}^{-1}y,y\rangle\varphi+\frac{\mathrm{tr}[K_{22,n}Q_{2,n}^{-1}]}{2}\varphi\\
&+\langle Q_{2,n}^{-1}y,K_{12,n}D_1\varphi\rangle.
\end{align*}
Remember that unitary transformation preserve essential m-dissipativity, therefore $(\hat{S}_n,\mathcal{S}(\R^n)\otimes C_0^{\infty}(\R^n))$ is essential m-dissipative in $L^2(\lambda^n\otimes\lambda^n)$.
\end{remark}
For the rest of this section we go back to the case, where our functions are \textbf{real valued}.\\
Our perturbation process is not done yet, indeed in the upcoming lines we study a perturbation with $(B_{\Psi},\mathcal{S}(\R^n)\otimes C_0^{\infty}(\R^n))$ defined by
\begin{align*}
\mathcal{S}(\R^n)\otimes C_0^{\infty}(\R^n)\ni \varphi\mapsto B_{\Psi}\varphi=\langle D \Psi(x),K_{21,n}D_2\varphi\rangle\in L^2(\lambda^n\otimes\lambda^n),
\end{align*}where $\Psi$ is Lipschitz continuous with global Lipschitz constant $d\in (0,\infty)$. Our goal is to set $\Psi(x)=\frac{1}{2}\langle Q_{1,n}^{-1}x,x\rangle$ for all $x\in\R$. But note that this particular $\Psi$ is only locally Lipschitz continuous. 
\begin{lemma}Suppose $\Psi$ is Lipschitz continuous with global Lipschitz-constant $d\in (0,\infty)$. The operator $(\hat{S}_n+B_{\Psi},\mathcal{S}(\R^n)\otimes C_0^{\infty}(\R^n))$ is essential m-dissipative in $L^2(\lambda^n\otimes\lambda^n)$.

\begin{proof}
First of all note that due to the standard integration by parts formula, the operator $(B_{\Psi},\mathcal{S}(\R^n)\otimes C_0^{\infty}(\R^n))$ is antisymmetric, hence dissipative with respect to the inner product of $L^2(\lambda^n\otimes\lambda^n)$.

Since $\Psi$ is Lipschitz continuous with global Lipschitz constant $d$ it holds
$s_{\Psi}=\norm{D\Psi}_{L^{\infty}(\lambda^n\otimes\lambda^n)}\leq \sqrt{d}n$. Set $b=\abs{K_{21,n}}$, $c=\mathrm{tr}[K_{22,n}Q_{2,n}^{-1}]$ and denote by $a_n$ the smallest eigenvalue of $K_{22,n}$.
We have for all $\varphi\in \mathcal{S}(\R^n)\otimes C_0^{\infty}(\R^n)$ using the antisymmetry of $(Q_{2,n}^{-1}y,K_{12,n}D_1)$
\begin{multline*}
\norm{B_{\Psi}\varphi}_{L^2(\lambda^n\otimes\lambda^n)}^2\leq \int_{\R^{n}\times\R^n}\abs{ D\Psi(x)}^2\abs{K_{21,n}}^2\abs{D_2\varphi}^2\mathrm{d}\lambda^n\otimes\lambda^n\\
\leq \frac{s_{\Psi}^2b^2}{a_n}\int_{\R^{n}\times\R^n} \abs{K_{22,n}^{\frac{1}{2}}D_2\varphi}^2\mathrm{d}\lambda^n\otimes\lambda^n= \frac{s_{\Psi}^2b^2}{a_n} (-\mathrm{tr}[K_{22,n}D_2^2\varphi],\varphi)_{L^2(\lambda^n\otimes\lambda^n)}\\
\leq \frac{s_{\Psi}^2b^2}{a_n}\big((- \hat{S}_n\varphi,\varphi)_{L^2(\lambda^n\otimes\lambda^n)}+\frac{c}{2}\norm{\varphi}^2_{L^2(\lambda^n\otimes\lambda^n)}\big).
\end{multline*}We can finish the proof by using \cite[Lemma 1.4.3.]{DissConrad}.
\end{proof}
\end{lemma}
For the next lemma we only assume that $\Psi$ is locally Lipschitz continuous. We use the strategy from \cite[Theorem 6.3.1.]{DissConrad} to show that $(\hat{S}_n+B_{\Psi},\mathcal{S}(\R^n)\otimes C_0^{\infty}(\R^n))$ is essential m-dissipative in $L^2(\lambda^n\otimes\lambda^n)$. The key in this consideration is a typical localization argument, therefore the following lemma introducing a collection of cut-off function is necessary.

\begin{lemma}\label{cutoff}Let $m\in\N$. There is some $\varphi\in C_0^{\infty}(\R^n)$ such that $0\leq\varphi \leq 1$, $\varphi=1$ on $B_1(0)=\lbrace x\in \R^n\vert\; \abs{x}< 1\rbrace$ and $\varphi= 0$ outside $B_2(0)$ and a constant $C\in(0,\infty)$, independent of $m\in \N$, such that 
\begin{align*}
\abs{\partial_i\varphi_m(z)}\leq\frac{C}{m},\quad \abs{\partial_{ij}\varphi_m(z)}\leq\frac{C}{m^2}\quad\text{for\;all}\quad z\in \R^n,\quad 1\leq i,j\leq n,
\end{align*}where we define
\begin{align*}
\varphi_m(z)=\varphi(\frac{z}{m})\quad\text{for\;each}\quad z\in \R^n,\quad m\in\N.
\end{align*}
In particular, $0\leq\varphi_m \leq 1$ and $\varphi_m=1$ on $B_m(0)$ for all $m\in \N$, besides $\varphi_m\rightarrow 1$ pointwisely on $\R^n$ and $D \varphi_m\rightarrow 0$ as $m\rightarrow\infty$, w.r.t. the sup norm $\abs{\cdot}_{\infty}$.
\end{lemma}

\begin{proposition}Assume that $\Psi$ is locally Lipschitz-continuous and $D\Psi\in L^2(\lambda^n)$, then $(\hat{S}_n+B_{\Psi},\mathcal{S}(\R^n)\otimes C_0^{\infty}(\R^n))$ is essential m-dissipative in $L^2(\lambda^n\otimes\lambda^n)$.
\begin{proof}
The dissipativity of the operator $(\hat{S}_n+B_{\Psi},\mathcal{S}(\R^n)\otimes C_0^{\infty}(\R^n))$ is derived by using the classical integration by parts formula. So it is left to show that $(Id-(\hat{S}_n+B_{\Psi}))(\mathcal{S}(\R^n)\otimes C_0^{\infty}(\R^n))$ is dense in $L^2(\lambda^n\otimes\lambda^n)$.

By the density of $C_0^{\infty}(\R^{n}\times\R^n)$ in $L^2(\lambda^n\otimes\lambda^n)$, it is enough to show that $C_0^{\infty}(\R^{n}\times\R^n)$ is contained in the closure of $(Id-(\hat{S}_n+B_{\Psi}))(\mathcal{S}(\R^n)\otimes C_0^{\infty}(\R^n))$.\\ Take $\varepsilon>0$ and an arbitrary $\psi\in C_0^{\infty}(\R^{n}\times\R^n)$. Now choose $\hat{\chi},\hat{\eta}\in C_0^{\infty}(\R^n)$ and define $\chi(x,y)=\hat{\chi}(x)$ and $\eta(x,y)=\hat{\eta}(x)$, for all $x,y\in\R^n$. Furthermore assume $0\leq\chi\leq\eta\leq 1$, $\chi=1$ on the support of $\psi$ and $\eta=1$ on the support of $\chi$.

Given $\varphi\in \mathcal{S}(\R^n)\otimes C_0^{\infty}(\R^n)$. It holds $D_1\chi \varphi=\varphi D_1\chi+\chi D_1\varphi$, combining this with the triangle inequality we get 
\begin{align*}
&\norm{(Id-(\hat{S}_n+B_{\Psi}))(\chi \varphi)-\psi}_{L^2(\lambda^n\otimes\lambda^n)}\\
=&\norm{\chi \big((Id-(\hat{S}_n+B_{\Psi}))\varphi-\psi\big)-\varphi\langle Q_{2,n}^{-1}y,K_{12,n}D_1\chi \rangle}_{L^2(\lambda^n\otimes\lambda^n)}\\
\leq & \norm{\chi \big((Id-(\hat{S}_n+B_{\Psi}))\varphi-\psi\big)}_{L^2(\lambda^n\otimes\lambda^n)}+\norm{\varphi\langle Q_{2,n}^{-1}y,K_{12,n}D_1\chi \rangle}_{L^2(\lambda^n\otimes\lambda^n)}.
\end{align*}By the support properties of $\chi,\eta$ and $\psi$ we have
\begin{align*}
 &\norm{\chi \big((Id-(\hat{S}_n+B_{\Psi}))\varphi-\psi\big)}_{L^2(\lambda^n\otimes\lambda^n)}\\
 =& \norm{\chi \big((Id-(\hat{S}_n+B_{\eta\Psi})\varphi-\psi\big)}_{L^2(\lambda^n\otimes\lambda^n)}\\
 \leq& \norm{(Id-(\hat{S}_n+B_{\eta\Psi}))\varphi-\psi}_{L^2(\lambda^n\otimes\lambda^n)}.
\end{align*}Since $\eta\Psi$ is Lipschitz-continuous we find $\varphi\in \mathcal{S}(\R^n)\otimes C_0^{\infty}(\R^n)$ such that 
\begin{align*}
\norm{(Id-(\hat{S}_{n}+B_{\eta\Psi})\varphi-\psi}_{L^2(\lambda^n\otimes\lambda^n)}< \frac{\varepsilon}{2}.
\end{align*}
Setting $a=\abs{K_{22,n}^{-\frac{1}{2}}}\abs{K_{12,n}}\abs{D_1\chi}_{\infty}$ and $c=\mathrm{tr}[K_{22,n}Q_{2,n}^{-1}]$ it holds
\begin{align*}
\norm{\varphi\langle Q_{2,n}^{-1}y,K_{12,n}D_1\chi \rangle}_{L^2(\lambda^n\otimes\lambda^n)}\leq a\norm{\abs{K_{22,n}^{\frac{1}{2}}Q_{2,n}^{-1}y}\varphi}_{L^2(\lambda^n\otimes\lambda^n)}
\end{align*}and for $b=\max\lbrace 4,2c\rbrace$
\begin{align*}
&\norm{\abs{K_{22,n}^{\frac{1}{2}}Q_{2,n}^{-1}y}\varphi}_{L^2(\lambda^n\otimes\lambda^n)}^2\\
\leq& 4(-(\hat{S}_{n}+B_{\eta\Psi})\varphi,\varphi)_{L^2(\lambda^n\otimes\lambda^n)}+4\frac{c}{2}(\varphi,\varphi)_{L^2(\lambda^n\otimes\lambda^n)}\\
\leq& b((Id-(\hat{S}_{n}+B_{\eta\Psi}))\varphi,\varphi)_{L^2(\lambda^n\otimes\lambda^n)}\leq b \norm{(Id-(\hat{S}_{n}+B_{\eta\Psi}))\varphi}_{L^2(\lambda^n\otimes\lambda^n)}^2,
\end{align*}where we used a similar argumentation as in Lemma \ref{Diss_fourier_pert} in the first inequality and the Cauchy-Schwarz inequality and the dissipativity of $\hat{S}_{n}+B_{\eta\Psi}$ in the third. 
A combination of the estimates above yields
\begin{align*}
&\norm{(Id-(\hat{S}_{n}+B_{\Psi})(\chi \varphi)-\psi}_{L^2(\lambda^n\otimes\lambda^n)}\\
\leq& \norm{(Id-(\hat{S}_{n}+B_{\eta\Psi})\varphi-\psi}_{L^2(\lambda^n\otimes\lambda^n)}+ a\norm{\abs{K_{22,n}^{\frac{1}{2}}Q_{2,n}^{-1}y}\varphi}_{L^2(\lambda^n\otimes\lambda^n)}\\
< &\frac{\varepsilon}{2}+ a\sqrt{b}\norm{ (Id-(\hat{S}_{n}+B_{\eta\Psi})\varphi}_{L^2(\lambda^n\otimes\lambda^n)}\\
\leq & \frac{\varepsilon}{2}+a\sqrt{b}\big(\norm{ (Id-(\hat{S}_{n}+B_{\eta\Psi}))\varphi-\psi}_{L^2(\lambda^n\otimes\lambda^n)}+\norm{\psi}\big)_{L^2(\lambda^n\otimes\lambda^n)}\\
\leq & \frac{\varepsilon}{2} + a\sqrt{b} \big( \frac{\varepsilon}{2}+\norm{\psi}_{L^2(\lambda^n\otimes\lambda^n)}\big).
\end{align*}Since by Lemma \ref{cutoff} we can choose $\chi$ such that \begin{align*}
a\sqrt{b} \big( \frac{\varepsilon}{2}+\norm{\psi}_{L^2(\lambda^n\otimes\lambda^n)}\big)<\frac{\varepsilon}{2},
\end{align*}we can conclude the assertion.
\end{proof}
\end{proposition}
Before we move on we fix
\begin{align*}
\Psi(x)=\frac{1}{2}\langle Q_{1,n}^{-1}x,x\rangle,\quad x\in\R^n,
\end{align*}to obtain essential m-dissipativity of $(\hat{S}_n+B_{\Psi},\mathcal{S}(\R^n)\otimes C_0^{\infty}(\R^n))$ in $L^2(\lambda^n\otimes\lambda^n)$. Note that for $\varphi\in \mathcal{S}(\R^n)\otimes C_0^{\infty}(\R^n)$ we have
\begin{align*}
(\hat{S}_n+B_{\Psi})\varphi=&\mathrm{tr}[K_{22,n}D_2^2\varphi]-\frac{1}{4}\langle Q_{2,n}^{-1}K_{22,n}Q_{2,n}^{-1}y,y\rangle \varphi+\frac{\mathrm{tr}[K_{22,n}Q_{2,n}^{-1}]}{2}\varphi\\
&-\langle Q_{1,n}^{-1}x,K_{21,n}D_2\varphi\rangle+ \langle Q_{2,n}^{-1}y,K_{12,n}D_1\varphi\rangle.
\end{align*}
\begin{theorem}\label{mdissipative}The operator $(L,\fc B_W))$ defined by
\begin{align*}
\fc B)\ni f\mapsto Lf=&\mathrm{tr}[K_{22}D^2_2f]-(v,Q_2^{-1}K_{22}D_2f)_V\\
&-(u,Q_1^{-1}K_{21}D_2f)_U+(v,Q_2^{-1}K_{12}D_1f)_V\in L^2(\mu),
\end{align*}is essentially m-dissipative in $L^2(\mu)$ with m-dissipative closure $(L,D(L))$. The resolvent in $\alpha\in (0,\infty)$, of $(L,D(L))$ is denoted by $R(\alpha,L)$.
\begin{proof}
Using the collection of cut-off functions provided by Lemma \ref{cutoff} and the Theorem of dominated convergence one can check that $C_0^{\infty}(\R^n)\otimes C_0^{\infty}(\R^n)$ is dense in $\mathcal{S}(\R^n)\otimes C_0^{\infty}(\R^n)$ w.r.t. the $\hat{S}_{n}+B_{\Psi}$-graph norm. Hence we obtain the essential m-dissipativity of $(\hat{S}_{n}+B_{\Psi},C_0^{\infty}(\R^n)\otimes C_0^{\infty}(\R^n))$ in $L^2(\lambda^n\otimes\lambda^n)$.\\
Recalling the unitary transformation $U:L^2(\lambda^n\otimes\lambda^n)\rightarrow L^2(\mu^n)$ one can calculate that for all $\varphi\in C_0^{\infty}(\R^n)\otimes C_0^{\infty}(\R^n))=U(C_0^{\infty}(\R^n)\otimes C_0^{\infty}(\R^n))$ we have
\begin{align*}
U(\hat{S}_{n}+B_{\Psi})U^{-1}\varphi=L_n\varphi.
\end{align*}Using the fact that unitary transformation preserves essential m-dissipativity we obtain the essential m-dissipativity of $(L_n,C_0^{\infty}(\R^n)\otimes C_0^{\infty}(\R^n))$ in $L^2(\mu^n)$. As $(L_n,C_b^{\infty}(\R^n\times\R^n))$ is a dissipative extension of $(L_n,C_0^{\infty}(\R^n)\otimes C_0^{\infty}(\R^n))$, we get essential m-dissipativity of $(L_n,C_b^{\infty}(\R^n\times\R^n))$ in $L^2(\mu^n)$.
Since by Lemma \ref{propL} we already know that $(L,\fc B_W))$ is dissipative it is left to show that $(Id-L)(\fc B_W))$ is dense in $L^2(\mu)$.

Let $\varepsilon>0$ and $h\in L^2(\mu)$ be given. Remark \ref{densitytwo} provides $g=\psi(P_n(\cdot),P_n(\cdot))\in \fc B_W,n)\subset L^2(\mu)$ such that
\begin{align*}
\norm{h-g}_{L^2(\mu)}< \frac{\varepsilon}{2}.
\end{align*}Since $(Id-L_n)(C_b^{\infty}(\R^n\times\R^n))$ is dense in $L^2(\mu^n)$ we find a $\varphi\in C_b^{\infty}(\R^n\times\R^n)$ such that 
\begin{align*}
\norm{(Id-L_n)\varphi-\psi}_{L^2(\mu^n)}< \frac{\varepsilon}{2}.
\end{align*}Set $f=\varphi(P_n(\cdot),P_n(\cdot))$ and use the triangle inequality together with Lemma \ref{imagemeasure} to observe
\begin{multline*}
\norm{(Id-L)f-h}_{L^2(\mu)}\leq \norm{(Id-L)f-g}_{L^2(\mu)}+\norm{h-g}_{L^2(\mu)}\\
=\norm{(Id-L_n)\varphi-\psi}_{L^2(\mu^n)}+\norm{h-g}_{L^2(\mu)}< \frac{\varepsilon}{2}+ \frac{\varepsilon}{2}=\varepsilon.
\end{multline*}
\end{proof}
\end{theorem}
Before we proof the main result of this section, we show a regularity estimate similar to the one from Theorem \ref{regestimatespotential}. In contrast to \ref{regestimatespotential} we don't have to deal with unbounded diffusion operators $(C,D(C))$ as coefficients, but the degenerate structure of $(L,\fc B_W))$ makes the proof more challenging.
\begin{proposition}For $f\in \fc B_W)$ and $\alpha \in (0,\infty)$, set $g=\alpha f-Lf$. Then the following equations hold
\begin{align*}
&\int_W\alpha f^2+\norm{K_{22}^{\frac{1}{2}}D_2f}^2_V\dm=\int_W fg\dm\\
&\int_W\alpha \norm{K_{22}^{\frac{1}{2}}D_2f}^2_V+\norm{Q_1^{-\frac{1}{2}}K_{21}D_2f}^2_U\\
&+\norm{Q_2^{-\frac{1}{2}}(K_{22}D_2f-K_{12}D_1f)}^2_V+\mathrm{tr}[(K_{22}D_2^2f)^2]\dm\\
&=\int_W(K_{22}D_2f,D_2g)_V-(K_{12}D_1f,D_2g)_V+(K_{21}D_2f,D_1g)_V\dm.
\end{align*}In particular we have
\begin{align*}
&\int_W\norm{Q_1^{-\frac{1}{2}}K_{21}D_2f}^2_U+\norm{Q_2^{-\frac{1}{2}}(K_{22}D_2f-K_{12}D_1f)}^2_V+\mathrm{tr}[(K_{22}D_2^2f)^2]\dm\\
&=\int_W (Lf)^2\dm.
\end{align*}
\begin{proof}
To show the first equation, we multiply $g=\alpha f-Lf$ with $f$, integrate over $W$ w.r.t. $\mu$ and use Lemma \ref{propL} item $(iv)$.\\
To derive the second equality we are taking partial derivatives in the direction of $d_i$ and $e_i$, $i\in\N$, separately to get
\begin{align*}
\alpha \partial_{i,1}f-L\partial_{i,1}f+(d_i,Q_1^{-1}K_{21}D_2f)_U&=\partial_{i,1}g\quad\text{and}\\
\alpha \partial_{i,2}f-L\partial_{i,2}f+(e_i,Q_2^{-1}K_{22}D_2f)_U-(e_i,Q_2^{-1}K_{12}D_1f)_U&=\partial_{i,2}g.
\end{align*}Multiply the first equation with $\partial_{j,2}f(d_i,K_{21}e_j)_U$ and sum over all $i,j\in\N$ to obtain
\begin{align*}
&\alpha (D_1f,K_{21}D_2f)_U-\sum_{i,j=1}^{\infty}L\partial_{i,1}f\partial_{j,2}f(d_i,K_{21}e_j)_U+(K_{21}D_2f,Q_1^{-1}K_{21}D_2f)_U\\
&=(K_{21}D_2f,D_1g)_V.
\end{align*}Similarly we multiply the second equation with $\partial_{j,2}f(e_i,K_{22}e_j)_V$ and sum over all $i,j\in\N$ yielding
\begin{multline*}
\alpha (D_2f,K_{22}D_2f)_V-\sum_{i,j=1}^{\infty}L\partial_{i,2}f\partial_{j,2}f(e_i,K_{22}e_j)_V\\
+(K_{22}D_2f,Q_2^{-1}K_{22}D_2f)_U-(K_{22}D_2f,Q_2^{-1}K_{12}D_1f)_V=(K_{22}D_2f,D_2g)_V.
\end{multline*}
Multiplying the second equation with $-\partial_{j,1}f(e_i,K_{12}d_j)_V$ and summing over all $i,j\in\N$ gives
\begin{multline*}
-\alpha (D_2f,K_{12}D_1f)_V+\sum_{i,j=1}^{\infty}L\partial_{i,2}f\partial_{j,1}f(e_i,K_{12}d_j)_U\\
-(K_{12}D_1f,Q_2^{-1}K_{22}D_2f)_V+(K_{12}D_1f,Q_2^{-1}K_{12}D_1f)_V=-(K_{12}D_1f,D_2g)_V.
\end{multline*}Summing up the equations we just derived and reordering the terms results in
\begin{multline*}
\alpha \norm{K_{22}^{\frac{1}{2}}D_2f}^2_V+\norm{Q_1^{-\frac{1}{2}}K_{21}D_2f}^2_U+\norm{Q_2^{-\frac{1}{2}}(K_{22}D_2f-K_{12}D_1f)}^2_V\\
+\sum_{i,j=1}^{\infty}L\partial_{i,2}f\partial_{j,1}f(e_i,K_{12}d_j)_V-L\partial_{i,1}f\partial_{j,2}f(d_i,K_{21}e_j)_U\\
-\sum_{i,j=1}^{\infty}L\partial_{i,2}f\partial_{j,2}f(e_i,K_{22}e_j)_V\\=(K_{22}D_2f,D_2g)_V-(K_{12}D_1f,D_2g)_V+(K_{21}D_2f,D_1g)_V.
\end{multline*}Integrating over $W$ w.r.t. $\mu$ together with the fact that
\begin{align*}
&\int_W\sum_{i,j=1}^{\infty}L\partial_{i,2}f\partial_{j,1}f(e_i,K_{12}d_j)_V-L\partial_{i,1}f\partial_{j,2}f(d_i,K_{21}e_j)_U\dm=0\quad\text{and}\\
&\int_W-\sum_{i,j=1}^{\infty}L\partial_{i,2}f\partial_{j,2}f(e_i,K_{22}e_j)_V\dm=\int_W\mathrm{tr}[(K_{22}D_2^2f)^2]\dm,
\end{align*}yields the second equation. As in Theorem \ref{regestimates} and Theorem \ref{regestimatespotential} we rearrange the terms in the second equation to get
\begin{align*}
&\int_W(K_{21}D_2f,Q_1^{-1}K_{21}D_2f)_U\\
&+(Q_2^{-1}(K_{22}D_2f-K_{12}D_1f),(K_{22}D_2f-K_{12}D_1f))_V+\mathrm{tr}[(K_{22}D_2^2f)^2]\dm\\
=&\int_W(K_{22}D_2f,D_2(-Lf))_V-(K_{12}D_1f,D_2(\alpha f-Lf))_V\\
&+(K_{21}D_2f,D_1(\alpha f-Lf))_V\dm\\
=&\int_W(K_{22}D_2f,D_2(-Lf))_V-(K_{12}D_1f,D_2(-Lf))_V\\
&+(K_{21}D_2f,D_1(-Lf))_V\dm=\int_W (Lf)^2\dm,
\end{align*}where the last equality is due to the integration by parts formula from Remark \ref{IBP_formel_ad} applied in $W^{1,2}(U,\mu_1,\R)$ and $W^{1,2}(V,\mu_2,\R)$.
\end{proof}
\end{proposition}
Proceeding as in Corollary \ref{domaindl} we obtain the following corollary.
\begin{corollary}For all $f\in D(L)$ it holds
\begin{align}\label{neuman}
\begin{split}
&\int_W\frac{1}{\nu_1}\norm{K_{21}D_2f}^2_U+\mathrm{tr}[(K_{22}D_2^2f)^2]\dm\\
&\leq\int_W\norm{Q_1^{-\frac{1}{2}}K_{21}D_2f}^2_U+\mathrm{tr}[(K_{22}D_2^2f)^2]\dm\leq \int_W(Lf)^2\dm.
\end{split}
\end{align}
\end{corollary}
\begin{theorem}\label{ess_diss_Lphi}Suppose $\norm{D\Phi}^2_{L^{\infty}(\mu_1)}< \frac{1}{4\nu_1}$, then the infinite-dimensional Langevin operator $(L_{\Phi},\fc B_W))$ is essential m-dissipative in $L^2(\mu^{\Phi})$.
\begin{proof}
By Lemma \ref{propL} item (iv) we already know that $(L_{\Phi},\fc B_W))$ is dissipative. In view of the Lumer-Phillips theorem we are left to show the dense range condition.
For $f\in L^2(\mu) $ set
\begin{align*}
Tf=-(D\Phi,K_{21}D_2R(1,L)f)_U.
\end{align*}We calculate using the Cauchy-Schwarz inequality, the assumption on $\Phi$ and Inequality \eqref{neuman}
\begin{align*}
\int_W (Tf)^2\dm &=\int_W (D\Phi,K_{21}D_2R(1,L)f)_U^2\dm\\
&\leq \norm{D\Phi}^2_{L^{\infty}(\mu_1)}\int_W (K_{21}D_2R(1,L)f,K_{21}D_2R(1,L)f)_U\dm\\
&\leq \norm{D\Phi}^2_{L^{\infty}(\mu_1)}\nu_1\int_W (Q_1^{-1}K_{21}D_2R(1,L)f,K_{21}D_2R(1,L)f)_U\dm\\
&<\frac{1}{4}\int_W (LR(1,L)f)^2\dm=\frac{1}{4}\int_W (f-R(1,L)f)^2\dm\leq \int_W f^2\dm.
\end{align*}Hence the operator $T:L^2(\mu)\rightarrow L^2(\mu)$ is well-defined. Moreover by Neumann's theorem we have $(Id-T)^{-1}\in \mathcal{L}(L^2(\mu))$. In particular for a given $g\in L^2(\mu)$ we find $f\in L^2(\mu)$ with $f-Tf=g$ in $L^2(\mu)$. Furthermore there is $h\in D(L)$ with $(Id-L)h=f$. Hence
\begin{align*}
(Id-L)h+(D\Phi,K_{21}D_2h)_U=f+(D\Phi,K_{21}D_2R(1,L)f)_U=f-Tf=g.
\end{align*}
This yields $L^2(\mu)\subset(Id-L_{\Phi})(D(L))$. Now similar as in Theorem \ref{essdispot} we can show that the closure  of $(Id-L_{\Phi})(\fc B_W))$ in $L^2(\mu^{\Phi})$ contains $(Id-L_{\Phi})(D(L))$. Since $L^2(\mu)$ is dense in $L^2(\mu^{\Phi})$ the dense range condition is shown.

\end{proof}
\end{theorem}

\section{Examples and Outlook}
In this section we have a look at certain examples, where the results we derived above are applicable. We consider the following situation, which is inspired by the one in \cite[Section 5]{SobolevRegularity}.\\
Let $U=V=L^2((0,1),\lambda,\R)$, $W=U\times V$, $K_{12}=K_{21}=K_{22}=Id$ and
$(-\Delta, D(\Delta))$ be the negative Dirichlet Laplacian, i.e. 
\begin{align*}
D(\Delta)&=W^{1,2}_0((0,1),\lambda,\R)\cap W^{2,2}((0,1),\lambda,\R)\subset L^2((0,1),\lambda,\R),\\
-\Delta x&=-x''.
\end{align*}
On $(U,\mathcal{B}(U))$ and $(V,\mathcal{B}(V))$ we consider two centered non-degenerate infinite-dimensional Gaussian measures $\mu_1$ and $\mu_2$ with covariance operators
\begin{align*}
Q_1=Q_2=-\Delta^{-1}:L^2((0,1),\lambda,\R)\rightarrow D(\Delta),
\end{align*}respectively. Recall the definition of $\mu_1^{\Phi}$ and $\mu^{\Phi}$ and denote by $B_U=B_V=(d_k)_{k\in\N}=(e_k)_{k\in\N}=(\sqrt{2}\sin(k\pi\cdot))_{k\in\N}$ the orthonormal basis of $L^2((0,1),\lambda,\R)$ diagonalizing $Q_1$ and $Q_2$ with corresponding eigenvalues $(\lambda_k)_{k\in\N}=(\nu_k)_{k\in\N}=(\frac{1}{k^2\pi^2})_{k\in\N}$. Additionally we fix a continuous differentiable (convex) function $\phi:\R\rightarrow \R$, which is bounded from below. Assume that there are constants $C_{1},C_2\in [0,\infty)$, $p_1\in [2,\infty)$ and $p_2\in [1,\infty)$ such that
\begin{align*}
&\abs{\phi(t)}\leq C_{1}(1+\abs{t}^{p_1}),\quad t\in\R,\\
&\abs{\phi'(t)}\leq C_{2}(1+\abs{t}^{p_2}),\quad t\in\R.
\end{align*}I.e.~$\phi$ and its derivative have at most polynomial growth. For such $\phi$ we consider potentials $\Phi:L^2((0,1),\lambda,\R)\rightarrow (-\infty,\infty]$ defined by
\begin{align*}
\Phi(u)=\begin{cases}\int_{(0,1)}\phi\circ u\mathrm{d}\lambda& u\in L^{p_1}((0,1),\lambda,\R) \\ \infty & u\notin L^{p_1}((0,1),\lambda,\R) \end{cases}.
\end{align*}
\begin{remark}\label{der}
Note that potentials as defined above are lower semicontinuous by Fatou's lemma, bounded from below and in $L^p(U,\mu_1,\R)$ for all $p\in [1,\infty)$. If $\phi$ is convex the same holds true for $\Phi$. Using \cite[Proposition 5.2]{SobolevRegularity} we know that $\Phi$ is bounded from below, lower semicontinuous and in $W^{1,2}(U,\mu_1,\R)$ with $D\Phi(u)=\phi'\circ u$ for a.e. $u\in L^2((0,1),\lambda,\R)$ (namely, for all $u\in L^{2p_2}((0,1),\lambda,\R)$).
\end{remark}
Choose a continuous differentiable (convex) $\tilde{\phi}$, which is bounded from below and with bounded derivative and define the potential $\Phi$ in terms of $\phi=\frac{\pi}{(2+\delta)\norm{\tilde{\phi}'}_{\infty}}\tilde{\phi}$, where $\delta >0$.
In \textbf{Hypothesis \ref{stass}} all items except the fourth are obviously fulfilled. But note that also the missing item is valid, since $(\Delta, D(\Delta))$ generates a strongly continuous contraction semigroup in $L^2((0,1),\lambda,\R)$ and for $\tau_1=\pi^2$ it holds
\begin{align*}
(-(-\Delta v),v)_{L^2((0,1),\lambda,\R)}\leq -\tau_1(v,v)_{L^2((0,1),\lambda,\R)},\quad v\in D(\Delta).
\end{align*}
By construction and Remark \ref{der} we have 
\begin{align*}
\norm{D\Phi}^2_{L^{\infty}(\mu_1)}=\norm{\frac{\pi}{(2+\delta)\norm{\tilde{\phi}'}_{\infty}}\tilde{\phi}'}^2_{L^{\infty}(\mu_1)}=\frac{\pi^2}{(2+\delta)^2}<\frac{\pi^2}{4}.
\end{align*}
Therefore \textbf{Theorem \ref{ess_diss_Lphi}} is applicable and we obtain essential m-dissipativity in $L^2(L^{2}((0,1),\lambda,\R)\times L^{2}((0,1),\lambda,\R),\mu,\R)$ of $(L_{\Phi},\fc B_W))$, where for $f\in \fc B_W)$ we have
\begin{align*}
L_{\Phi}f=&\mathrm{tr}[D^2_2f]+(v,\Delta D_2f)_{L^2(\lambda)}\\
&+(u,\Delta D_2f)_{L^2(\lambda)}-(\phi'\circ u, D_2f)_{L^2(\lambda)}-(v,\Delta D_1f)_{L^2(\lambda)}.
\end{align*}
This is the starting point to construct a martingale and even a weak solution to the non-linear infinite-dimensional stochastic differential equation given by
\begin{align}\begin{split}\label{sdeex}
	&\mathrm{d}U_t=-\Delta V_t\mathrm{d}t\\
	&\mathrm{d}V_t=(\Delta V_t+\Delta U_t-\phi'(U_t))\mathrm{d}t+\sqrt{2}\mathrm{d}W_t.
	\end{split}
\end{align}
In the equation above $(W_t)_{t\geq 0}$ is a cylindrical Brownian motion on $(V,\mathcal{B}(V))$. 
A heuristically application of the Itô-formula suggest that this stochastic differential equation corresponds to $(L_{\Phi},D(L_{\Phi}))$. In order to make this correspondence rigorous and to construct weak solutions we plan to apply general resolvent methods described in \cite{boboc}.

\bigskip

As announced in the introduction the results derived in Section \ref{section_reg} are naturally needed while applying the abstract Hilbert space hypocoercivity method from \cite{GrothausHypo}. A rigorous application of this method in our infinite-dimensional setting goes beyond the aim of this article. At this point we just point out, how Ornstein-Uhlenbeck operators perturbed by the gradient of a potential with unbounded diffusions $(C,D(C))$ as coefficients appear during the application process. To this end define $P_S:L^2(\mu^{\Phi})\rightarrow L^2(\mu^{\Phi})$ by
\begin{align*}
P_Sf= \int _V f\mathrm{d}\mu_2,
\end{align*}where the integration is understood w.r.t. the second variable. An application of Fubini's theorem and the fact that $(V,\mathcal{B}(V),\mu_2)$ is a probability space shows that $P_S$ is a well-defined orthogonal projection on $L^2(\mu^{\Phi})$ satisfying
\begin{align*}
P_Sf\in L^2(\mu_1^{\Phi})\quad \text{and}\quad \norm{P_Sf}_{L^2(\mu_1^{\Phi})}= \norm{P_Sf}_{L^2(\mu^{\Phi})},\quad f\in L^2(\mu^{\Phi}).
\end{align*} 
In the definition above we canonically embed $L^2(\mu_1^{\Phi})$ into $L^2(\mu^{\Phi})$. Using that $\mu^{\Phi}$ is a probability measure one can check that the map $P:L^2(\mu^{\Phi})\rightarrow L^2(\mu^{\Phi})$ given as
\begin{align*}
Pf=P_Sf-(f,1)_{L^2(\mu^{\Phi})},\quad f\in L^2(\mu^{\Phi}),
\end{align*}is an orthogonal projection with
\begin{align*}
Pf\in L^2(\mu_1^{\Phi})\quad\text{and}\quad \norm{Pf}_{L^2(\mu_1^{\Phi})}= \norm{Pf}_{L^2(\mu^{\Phi})},\quad f\in L^2(\mu^{\Phi}).
\end{align*} One can show that for all $f\in \fc B_W)$, the operator $(PA_{\Phi}^2P,\fc B_W))$ is given by the formula
\begin{align*}
PA_{\Phi}^2Pf=&\mathrm{tr}[CD^2f_S]-(u,Q_1^{-1}CD f_S)_U-(D\Phi(u),CD f_S(u))_U,
\end{align*} where $f_S=P_Sf\in \fc B_U)$ and the possible unbounded operator $(C,D(C))$ in $U$ is given by
\begin{align*}
C=K_{21}Q_2^{-1}K_{12}\quad \text{with}\quad D(C)=\lbrace u\in U\vert K_{12}u\in D(Q_2^{-1})\rbrace.
\end{align*} 
The essential m-dissipativity of $(PA_{\Phi}^2P,\fc B_W))$ in $L^2(\mu^{\Phi})$ and corresponding regularity estimates, which are applicable in view of Theorem \ref{essdispot} and Theorem \ref{theorem_reg} derived in Section \ref{section_reg}, are fundamental to show \cite[Corollary 2.13]{GrothausHypo} and \cite[Proposition 2.15]{GrothausHypo}.

In a future article we would like to apply the general abstract hypocoercivity method in our infinite-dimensional setting. I.e.~we want to show that the semigroup $(T_t)_{t\geq 0}$ generated by $(L_{\Phi},D(L_{\Phi}))$ is hypocoercive (compare \cite[Theorem 2.18]{GrothausHypo}). Furthermore, using the methods from \cite{boboc} we plan to construct martingale/weak solutions to \eqref{sdeex}. Additionally our aim is to combine these results to study the longtime behavior of the martingale/weak solution by relating the transition semigroup corresponding to \eqref{sdeex} with $(T_t)_{t\geq 0}$.







\end{document}